\documentclass[12pt,a4paper,reqno]{amsart}

\usepackage{amsmath,amsthm,amsfonts,amssymb,euscript,cite, mathrsfs,mathtools}
\usepackage{graphics}
\usepackage{xcolor} %
\usepackage{hyperref}
\hypersetup{colorlinks,linkcolor={blue},citecolor={blue},urlcolor={red}}  
\usepackage{enumitem} %
\usepackage{ulem}
\usepackage[noabbrev,capitalize,nameinlink]{cleveref}
\usepackage{setspace}
\usepackage{enumitem}

\crefname{equation}{}{}
\crefname{algocf}{Algorithm}{Algorithms}
\crefname{equation}{}{} %
\crefname{algocf}{Algorithm}{Algorithms}

\usepackage{graphicx}
\usepackage{color}
\usepackage{transparent}

\usepackage{fullpage} %

\usepackage{microtype} %

\usepackage{seqsplit}

\definecolor{green}{rgb}{0,0.8,0} %

\definecolor{babypink}{rgb}{0.96,0.76,0.76}

\newcommand{\ang}{{\not\negmedspace\partial }}

\newcommand{\rt}{{\mathbb R^3}}

\newcommand{\pa}{\partial}

\newcommand{\grad}{\nabla_x}

\newcommand{\gab}{g^{\alpha\beta}}
\newcommand{\hab}{h^{\alpha\beta}}

\newcommand{\pab}{\partial_\beta}
\newcommand{\paa}{\partial_\alpha}

\newcommand{\pat}{\partial_t}
\newcommand{\pav}{\pa_v}

\newcommand{\pao}{\partial\mkern-10mu /\,}

\newcommand{\ti}{\tilde}

\newcommand{\la}{\langle}
\newcommand{\ra}{\rangle}

\newcommand{\ls}{\lesssim}
\newcommand{\gs}{\gtrsim}

\newcommand{\inv}{^{-1}}
\newcommand{\invh}{^{-\f12}}

\newcommand{\de}{\nabla_{t,x}} %

\newcommand{\f}{\frac}

\newcommand{\iy}{\infty}
\renewcommand{\S}{{\mathbb S}}

\newcommand{\ltst}{{L^2(S^2)}}
\newcommand{\crt}{{C^{ R}_{ T}}}
\newcommand{\cut}{{C^{ U}_{ T}}}
\newcommand{\ctr}{{C^T_R}}

\newcommand{\crtt}{{\ti C^{ R}_{ T}}}

\newcommand{\inte}{{C^{<3T/4}_{ T}}}

\newcommand{\co}{{D_{tr}^{R}}}
\newcommand{\dtr}{{D_{tr}}}

\newcommand{\lolt}{{L^1L^{2}}}
\newcommand{\ltlt}{{L^2L^2}}

\newcommand{\lt}{L^2}
\newcommand{\p}{\phi}

\newcommand{\supp}{\text{supp\,}}

\let\pminus\pm
\renewcommand{\pm}{\phi_{\le m}}
\newcommand{\pn}{\phi_{\le n}}

\newcommand{\pmn}{\p_{\le m+n}}
\newcommand{\nm}{\la t-r\ra} %
\newcommand{\jr}{\la r\ra} %
\newcommand{\ju}{\la u\ra}

\newcommand{\jt}{\la t\ra}
\newcommand{\jv}{\la v\ra}
\newcommand{\js}{\la s\ra}
\newcommand{\jrho}{\la\rho\ra}

\newcommand{\lr}[1]{\left( #1 \right)}

\newcommand{\tpsi}{{\widetilde{\psi}}}
\newcommand{\les}{{LE^*}}
\newcommand{\leo}{{LE^1}}

\def\doi#1{ {\href{http://dx.doi.org/#1}
   {{\mdseries\ttfamily DOI}}}}

\newtheorem{theorem}{Theorem}[section]
\newtheorem{corollary}[theorem]{Corollary}
\newtheorem{lemma}[theorem]{Lemma}
\newtheorem{proposition}[theorem]{Proposition}
\theoremstyle{definition}
\crefname{claim}{Claim}{Claims}

\newtheorem{definition}[theorem]{Definition}
\newtheorem{example}[theorem]{Example}

\theoremstyle{remark}
\newtheorem{remark}[theorem]{Remark}
\theoremstyle{conjecture}

\numberwithin{equation}{section}

\newcommand{\x}{\alpha}
\newcommand{\xb}{\beta}
\newcommand{\xd}{\delta}
	
\newcommand{\xg}{\gamma}
	
\newcommand{\eps}{\epsilon}

\newcommand{\xk}{\kappa}
\newcommand{\xl}{\lambda}
	
\newcommand{\xo}{\omega}
	\newcommand{\xO}{\Omega}

\newcommand{\xs}{\sigma}
	\newcommand{\xS}{\Sigma}

\newcommand{\N}{{\mathbb N}}

\newcommand{\R}{\mathbb R}

\newcommand{\calC}{\mathcal C}

\newcommand{\calN}{\mathcal N}

\newcommand{\calR}{\mathcal R}

\newcommand{\calT}{\mathcal T}

\begin{document}

\title{Decay rates for cubic and higher order nonlinear wave equations on asymptotically flat spacetimes} %
\author{Shi-Zhuo Looi}
\address{Department of Mathematics, University of Kentucky, Lexington, 
  KY  40506}

\begin{abstract}In this paper, we prove pointwise decay rates for cubic and higher order nonlinear wave equations, including quasilinear wave equations, on asymptotically flat and time-dependent spacetimes. We assume that the solution to the linear equation (rather than the nonlinear equation) satisfies a weaker form of the standard integrated local energy decay, or Morawetz, estimate. For nonlinearities with a total derivative structure, we prove better pointwise decay rates.
\end{abstract}

\maketitle

\section{Introduction}

In this paper we study the long-time behaviour of solutions to cubic and higher order nonlinear wave equations in three spatial dimensions on time-dependent and asymptotically flat spacetimes. That is, we study these wave equations on curved spacetimes. Such equations can be viewed as perturbations of the classical linear wave equation, at least for small solutions which also have sufficiently small derivatives. 
	The goal is to obtain certain pointwise decay rates stated in \cref{thm:main} and then to prove modified decay rates for certain exceptional nonlinearities stated in \cref{thm:main2}; this is achieved by an iteration scheme that is outlined in \cref{outline of iteration}. 

With the exception of the cubic, quartic and quintic power-type nonlinear equations, which were studied in \cite{Loo22}, this article proves pointwise decay rates for cubic and higher order nonlinear wave equations, on asymptotically flat spacetimes with minimal assumptions.
	We note that quadratic %
equations that contain derivatives were considered in the work \cite{LooToh22}, where analogous long-time behaviour questions related to pointwise decay were investigated.

We consider the operator \begin{equation}\label{P def}
P := \paa\gab(t,x)\pab + g^\xo(t,x) \Delta_\xo + B^\x(t,x)\paa + V(t,x) \quad \text{on }\R^{1+3}
\end{equation} describing our spacetimes, where the coefficients are allowed to depend on $t$ and we use the summation convention for repeated indices. 
Here $\Delta_\xo$ denotes the Laplace operator on the unit sphere, and $\x,\xb$ range across $0, \dots, 3$.  See \cref{P.assptns} for the precise assumptions on $P$.  
	We study the nonlinear Cauchy problem
\begin{equation}   \label{eq:problem}
\begin{cases}
  P\phi(t,x) =  \calC(\pa^{2}\p,\pa\p,\p) & (t,x) \in \R^{1+3}\\
(\p(0,x), \pat \p(0,x)) =(\p_0,\p_1)
\end{cases}
\end{equation}where $\calC$ is an $n$-th order form whose terms are all at least cubic; see \cref{def:C} for the precise definition of $\calC$. 
	Our main theorem (\cref{thm:main}) states, informally, that if the coefficients of $P - \Box$ are asymptotically flat and if an integrated local energy decay statement (i.e., a Morawetz estimate) is satisfied, then the solution to \eqref{eq:problem} (and its vector fields) obey the global pointwise decay rates of $|\p(t,x)| \le C\la t-r\ra^{-1 - \min(c(P), A)}\la t+r\ra^{-1}$. Here $c(P)$ is a constant depending on the coefficients in \cref{P def} while $A$ is a constant depending on the nonlinearity. Thus for bounded $|x|$, we have $|\p(t,x)| \le C t^{-2-\min(c(P),A)}.$ We allow for large and non-compactly supported initial data in a localised $\lt$ space: see \cref{data.assptns}.
	An overview of the proof is contained in \cref{ss:outline}.
	
\subsection*{History of related problems}
In \cite{KP}, a similar class of equations was considered on the flat spacetime, i.e. on Minkowski spacetime. The work \cite{KP} used elementary $L^p$-$L^q$ estimates in certain Sobolev spaces from \cite{Pec,vonW} and avoided the need to use a Nash-Moser-type iteration scheme which was used in the earlier works \cite{K1,K3}. The work \cite{Sha} also uses $L^p$-$L^q$ estimates to prove the same results as those obtained in \cite{KP}, and it uses a global Picard iteration scheme. 

For quadratic equations, the semilinear wave equation in $\R^{1+3}$
\begin{equation}
\Box \phi = Q(\pa\p,\pa\p), \qquad \phi |_{t=0} = \phi_0, \qquad \pa_t \phi |_{t=0} = \phi_1
\end{equation}
for small initial data was been studied extensively. It is known that the solution blows up in finite time if $Q(\pa\p,\pa\p) = (\pa_t\p)^2$. On the other hand, if the nonlinearity satisfies the null condition (see, for instance, \cite{LooToh22} for a definition) it was shown in \cite{Ch} and \cite{K2} that the solution exists globally. This result was extended to quasilinear systems with multiple speeds, as well as the case of exterior domains (see, for instance, \cite{MNS}, \cite{MS1}) and to systems satisfying the weak null condition, including Einstein's Equations (see \cite{LR2}, \cite{LR3}). See also the upcoming \cite{LukOh} for sharp pointwise bounds and asymptotics, given certain assumptions.

The theory of global existence, uniqueness and scattering for the semilinear wave equation in $\R^{1+3}$
$$\Box \phi = \pminus \p^{p+1}, \qquad \p(0,x) = \phi_0(x), \qquad \pa_t \phi(0,x) = \phi_1(x)$$ 
was studied extensively; for instance, in the articles \cite{BSh,Gri}.
	Work has also been done for the pointwise decay of solutions; see \cite{Pec1,Yang}. In the case of compactly supported smooth data, decay rates were proved in \cite{Sz} (for small data) and in \cite{Gri,BSz} (for large data).

\subsection*{Setup and statement of main theorem}
We state some notation that we use throughout the paper.
	We write $X\ls Y$ to denote $|X| \leq CY$ for an implicit constant $C$ which may vary by line. Similarly, $X \ll Y$ will denote $|X| \le c Y$ for a sufficiently small constant $c>0$. 
	\, In $\R^{1+3}$, we consider %
\[
\partial := (\partial_t, \partial_1,\pa_2,\pa_3), \qquad \Omega := (x^i \partial_j -
x^j \partial_i)_{i,j}, \qquad S := t \partial_t + \sum_{i=1}^3x^i \partial_{i},
\]
which are, respectively, the generators of translations, rotations and scaling. We denote the angular derivatives by $\ang$.
	We set
$$Z := (\pa,\Omega,S) $$ and we define the function class $$S^Z(f)$$ to be the collection of real-valued functions $g$ such that $|Z^J g(t,x)| \ls_J |f|$
whenever $J$ is a multiindex. We will frequently use $f = \jr^\x$ for some real $\x \leq 0$, where $\jr := (1 + |r|^2)^{1/2}$. We also define $S^Z_\text{radial}(f) := \{ g \in S^Z(f) : g \text{ is spherically symmetric}\}.$
	We denote
\begin{align}  \label{vf defn}
\begin{split}
\p_{J} &:= Z^J\p := \partial^i \Omega^j S^k u, \  \text{ if } J = (i,j,k)\\
\pm &:= (\p_{J})_{|J|\le m}, \quad \p_{=m} := (\p_{J})_{|J| = m}.
\end{split}
\end{align}Thus $\p_{\vec0} = \phi$.  %
We let
\begin{equation}
\label{new.def}S^{Z}_{\pa}(1) := \{ f \in S^{Z}(1) : \pa f \in S^{Z}(\jr\inv) \}.
\end{equation}For instance, constants belong in the class $S^Z_{\pa}(1)$. 
We let
\begin{equation}
\label{def:tangential}\bar\pa := \{ \pat+\pa_{r}, \pao\}
\end{equation}
where $\pao$ denotes the angular derivatives.

Forward energy bounds are usually viewed as a necessary prerequisite for local energy decay and pointwise bounds. We shall assume our equation \cref{eq:problem} satisfies the following weaker version of energy bounds, where we allow losses on the right hand side:
\begin{definition}[Weak energy bounds] \label{def:EB}
	 		We will assume that \eqref{eq:problem} satisfies the following estimate: there exists some $k_0 \in \N$ such that for finitely many $m \in\N$,
\begin{equation}\label{EB}
 \| \pa\p(T_{1}) \|_{H^{m}( \rt))} 
\ls_m  \|\pa\p(T_0)\|_{H^{m+k_{0}}(\rt)}, \quad 0 \le T_{0} \le T_{1}.
\end{equation}

\end{definition}%

\subsection*{Local energy norms}
The usual local energy decay estimate is as follows in \cref{LED}. Before we state the estimate, we define the $\leo$ norm. Let
$$A_R := \{ x\in\rt:R<|x|<2R\} \ \ (R>2), \qquad A_{R=1} := \{ |x|<2 \}.$$
Given a subinterval $I$ of $\R^+$,
\begin{equation}\label{initial.LE.def}
\begin{split}
 \| \phi\|_{LE(I)} &:= \sup_R  \| \la r\ra^{-\frac12} \phi\|_{L^2 (I\times A_R)},\\
  \| \phi\|_{LE^1(I)} &:= \| \de \phi\|_{LE(I)} + \| \la r\ra^{-1} \phi\|_{LE(I)},\\
 \| f\|_{LE^*(I)} &:= \sum_R  \| \la r\ra^{\frac12} f\|_{L^2 (I \times A_R)}.
\end{split} 
\end{equation}
Higher-order versions of \cref{initial.LE.def} are as follows:
\[
\begin{split}
  \| \phi\|_{LE^{1,k}(I)} &= \sum_{|\alpha| \leq k} \| \partial^\x \phi\|_{LE^1(I)} \\
  \| \phi\|_{LE^{0,k}(I)} &= \sum_{|\alpha| \leq k} \| \partial^\x \phi\|_{LE(I)},\\
  \| f\|_{LE^{*,k}(I)} &=  \sum_{|\alpha| \leq k}  \| \partial^\alpha f\|_{LE^{*}(I)}.
\end{split}  
\]
If the subinterval $I$ is omitted, then the norm will involve an integration over $\R^{+}$.

We have the following scale-invariant estimate on Minkowski backgrounds:
\begin{equation}\label{localenergyflat}
\|\pa \p\|_{L^{\infty}_t L^2_x} + \| \p\|_{LE^1}
 \lesssim \|\pa\p(0)\|_{L^2} + \|\Box \p\|_{LE^*+L^1_t L^2_x}
\end{equation}
and a similar estimate involving the $LE^1[t_0, t_1]$ and $LE^*[t_0, t_1]$ norms.
	See \cite{M} for the Klein-Gordon case, and see the following works in the case of small perturbations of the Minkowski space-time: see for instance \cite{KSS}, \cite{KPV}, \cite{SmSo}. Even for large perturbations, in the absence of trapping, \eqref{localenergyflat} still sometimes holds, see for instance \cite{BH}, \cite{MST}.
In the presence of trapping, \eqref{localenergyflat} is known to fail, see \cite{Ral}, \cite{Sb}. 
We will assume that a similar, but weaker, estimate holds (see \cref{def:SILED}), for our operator $P$ after commuting with \textit{only a finite number of derivatives} (as opposed to vector fields $Z$). Moreover, we do not assume that we can control the time derivative on the left-hand side---see \cref{eq:SILED}.

\begin{definition}[ILED] \label{LED}
We say that $P$ has the integrated local energy decay
property if the following estimate holds for a finite number $m \geq 0$, and for all $0\leq T_0<T_1 \leq \infty$:
\begin{equation}
 \| \pa^{\le m} \p\|_{LE^{1}([T_0,T_1)\times\rt)} 
\ls_m  \|\pa \p(T_0)\|_{H^m(\rt)} + \|\pa^{\le m}(P\p)\|_{(\lolt + LE^*)([T_0,T_1)\times\rt)}
\label{eq:LED}\end{equation}
where the implicit constant does not depend on $T_0$ and $T_1$.
\end{definition}

We now define a weaker version of \cref{LED}. This is the Morawetz estimate that we shall assume for the equation \cref{eq:problem}. 
   \begin{definition}[SILED] \label{def:SILED}
	 		We say that $P$ has the stationary integrated local energy decay property if the following estimate holds for a finite number of $m \ge 0$, and for all $0 \le T_0 < T_1 \le\iy$:
\begin{equation}\label{eq:SILED}
 \| \pa^{\le m}\p\|_{LE^1([T_0,T_1)\times\rt)} 
\ls_m  \|\pa \p(T_0)\|_{H^{m}(\rt)} + \|\pa^{\le m}(P\p)\|_{ LE^*([T_0,T_1)\times\rt)} + \|\pat \pa^{\le m} \p\|_{LE([T_0,T_1)\times\rt)}.
\end{equation}
\end{definition}Clearly the estimate \cref{eq:SILED} is a weaker assumption than the ILED estimate \cref{eq:LED}. The ILED estimate is a standard tool in the field, and has been proven using a variety of tools, and moreover in various different settings.

For the Kerr spacetime with small angular momentum, \cref{eq:LED,eq:SILED} have both been proven.
	For the Kerr spacetime with large angular momentum, the estimate \cref{eq:LED} holds for the homogeneous wave equation: see Theorem 3.2 in \cite{DRS}.
	The estimate \cref{eq:SILED} was proven in Theorem 4.3 of the article \cite{MTT}, for perturbations of the Schwarzschild spacetime that include the Kerr spacetime with small angular momentum. 
	In that result, no decay in time %
near the trapped set (which is the main issue) was even necessary for the metric.  In the stationary case, \cref{eq:SILED} can be thought of as a substitute for an elliptic estimate at zero frequency.

\begin{definition}[Definition of the nonlinearity $\calC$ being studied]\label{remark on C}\label{def:C}
Let $\calC(x,y,z)$ denote a polynomial of degree $\text{deg} \ge 3$ whose terms are all at least cubic, with $z^{5}, z^{4}, z^{3}$ disallowed. Thus $\calC(\pa^{2}\p, \pa\p,\p)$ allows all cubic and higher order nonlinearities except for terms of the form $\phi^3, \phi^{4}, \p^{5}$.\footnote{The cases of the nonlinearities $\p^{3},\p^{4},\p^{5}$ were covered in the article \cite{Loo22}. See remarks in the introduction and the conclusion of \cite{Loo22} for commentary on the cubic and quartic cases.} 
\begin{itemize}
\item
Let $\calN \in \N_{\ge3}$ denote the order of the nonlinearity $P\p$.%
\item
Let $\calT \in \N_{0}$ be the number of tangential derivatives that is present in \textit{each} term of the nonlinearity $P\p$. The tangential derivatives were defined in \cref{def:tangential}.
\end{itemize}
\end{definition}

\begin{example}[Examples of \cref{def:C}]
For instance, $(\calT,\calN)=(2,3)$ for $P\phi = \p \bar\pa \p \bar\pa \p$ (since there are two tangential derivatives), while for $$P\p = \p\pa\phi\pa\phi$$ %
 we have $(\calT,\calN) = (0,3)$. \  For another example, consider the classical null condition, for which $\calT=1$ (since each term in the null form has at least one tangential derivative) and $\calN=2$ (since each term is a product of two functions).  \ For $P\p = (\pa\p)^{3} + (\pa \p)^{4}$, we have $(\calT,\calN) = (0,3)$. 
\end{example}

\subsubsection{Assumptions on $P$}\label{P.assptns}
Let $h=g-m$, where $m$ denotes the Minkowski metric. Let $\xs \in (0,\infty)$ be real. We make the following assumptions on the coefficients of $P$: 
\begin{equation}\label{coeff.assu}
\begin{split}
 h^{\alpha\beta}, B^{\alpha} \in S^Z(\jr^{-1-\xs}) \\
\pat B^{\alpha}, V\in S^Z(\jr^{-2-\xs}) \\
g^\xo \in S^Z_\text{radial}(\jr^{-2-\xs});
\end{split}
\end{equation}that is, these are ``rough'' backgrounds which apply to a wide variety of situations.

In the special case of nonlinearities that contain the cubic form $c(t,x)\phi^{2}\pa\p$ with bounded coefficients $c(t,x) \in S^{Z}(1)$ (i.e. $c$ and all its vector fields are bounded), we assume that $\xs > 1/4$. More precisely, in this case we can make the assumption 
\begin{equation*}
 h^{\alpha\beta}, B^{\alpha} \in S^Z(\jr^{-1-\xs}), \ 
\pat B^{\alpha}\in S^Z(\jr^{-2-\xs}), \ 
V \in S^{Z}(\jr^{-2-\eps}), \ \eps >0 \text{ arbitrary}.
\end{equation*}

In what follows, we state our main theorem (\cref{thm:main}), as well as a follow-up theorem (\cref{thm:main2}).
\begin{theorem}\label{thm:main}Let $\phi$ solve \cref{eq:problem} with the assumptions \cref{coeff.assu} and let
$$ \kappa := \min(\sigma,\calT + \calN - 3).$$
Fix $m \in \N$. We assume $\p_0\in \lt(\rt)$ and that for a fixed $N \gg m$,  
\begin{equation}
\label{data.assptns}
	\sum_{J : |J|=0}^{N}\|\jr^{1/2 + \kappa}\pa\p_J(0)\|_{\lt(\rt)} + \|\p_{J}\|_{LE^{1}}<\iy
\end{equation}	
and that the solution exists globally (for instance, this holds in the small data case).
We assume \cref{EB,eq:SILED}. Then
we have the pointwise decay rate
\begin{equation}\label{eq:main bound}
\sum_{J:|J|=0}^m |\p_J(t,x)| \ls \f1{\la t+|x|\ra\la t-|x|\ra^{1+\kappa}}.
\end{equation}

\end{theorem}
	Thus given a fixed $x$, the solution $\p$ and its vector fields obey the pointwise upper bound $Ct^{-2-\kappa}$, and for solutions to \cref{eq:eqn,eq:eqn2}, we obtain faster decay than this rate by $t\inv$ if the coefficients of $P-\Box$ decay sufficiently rapidly (e.g., if on the Minkowski spacetime, where $P-\Box$ in fact equals zero). 

\begin{remark}\label{incl.quadr}
The pointwise decay rate \cref{eq:main bound} also holds for other nonlinearities, including quadratic derivative and power-type nonlinearities. 
For the power-type nonlinearities $P\p = \pminus \p^{p+1}$ or $P\p = |\p|^{p}\p$, with $p\in \N_{\ge2}$, the rate \cref{eq:main bound}, which reads $\la t+|x|\ra \inv \la t-|x|\ra^{- ( 1+\min(\xs, p - 2) )}$ since $(\calT,\calN)=(0,p+1)$, was proved in \cite{Loo22}.%

The decay rate \cref{eq:main bound} in fact also holds for both the semilinear and quasilinear wave equations with a quadratic nonlinearity satisfying the classical null condition. In that case, $\calN=2$ while $\calT=1$, thus \cref{eq:main bound} gives $\pm \ls \jv\inv\ju\inv$, which matches the sharp decay rate for the classical null condition found in \cite{LooToh22}. %

Note however that in the special case that the nonlinearity possesses a total derivative structure, then one can do better than the 
rate \cref{eq:main bound}. This is explained in \cref{thm:main2} and its proof (see \cref{sec:final}).
\end{remark}

\begin{theorem} \label{thm:main2}
\begin{enumerate}	
\item Let $P$ be defined as in \cref{P def} with spherically symmetric coefficients.

For spherically symmetric solutions of the semilinear wave equation\footnote{See \cref{new.def} for the definition of $S^{Z}_{\pa}(1)$.} 
\begin{equation}\label{eq:eqn}
P\phi (t,r)= \sum_{i=1}^{M} c_{i}(t,r) \phi^{n}\pa_{(i)}\p, \quad \pa_{(i)}\in \{\paa\}_{\x \in \{0,\dots,3\}}, 
\ \ c_{i}\in S^{Z}_{\pa}(1), \quad n \ge2, \ M \ge1
\end{equation} we have better decay than %
\cref{eq:main bound} by $\nm\inv$ for the (nonlinear component of this) equation:
\begin{equation}\label{eq:main bound2}
|\p(t,r)|\ls \f1{\la t+r\ra\la t-r\ra^{\min(1+\xs,\calT+\calN-1)}} = \f1{\la t+r\ra\nm^{\min(1+\xs,\calT+n)}}.
\end{equation}
For instance, for bounded $r$ we have
$$|\p(t,r)|\le C (1+t)^{-(\calT+n+1)} + C(1+t)^{-2-\sigma}.$$
For the $n=2$ case of the semilinear equation \cref{eq:eqn}, we assume that $\|\p_{\le N}\|_\leo$ is sufficiently small.%

\ \ In addition, the same decay rate \cref{eq:main bound2} also holds for spherically symmetric solutions of
the quasilinear wave equation
\begin{equation}\label{eq:eqn2}
P\phi (t,r)= \sum_{i=1}^{M} c_{i}(t,r) Y^{n}\pa_{(i)}Y, \quad Y \in \{ \pa \p \}
\ \ c_{i}\in S^{Z}_{\pa}(1), \quad n \ge2, \ M \ge1
\end{equation}%

\item Furthermore, for both \cref{eq:eqn,eq:eqn2}, 
if $\pa_{(i)} = \pat$, then the bound \cref{eq:main bound2} holds for solutions $\p$, coefficients $c_{i}$ and operators $P$ that are not necessarily spherically symmetric, thus $\p = \p(t,x)$, $c_{i} = c_{i}(t,x)$ and $(P-\Box) = (P-\Box)(t,x)$. 
\end{enumerate}

\end{theorem}

\begin{remark}
The argument shown in this paper straightforwardly yields a proof of a more general version of \cref{thm:main,thm:main2} when the decay increments $\sigma$ differ for the coefficients. (In \cref{coeff.assu}, the increments are all assumed to be equal to $\sigma$.) See also the main theorem in \cite{L} for a demonstration of this claim.
\end{remark}

\begin{remark}
All the arguments in this paper can be adapted to the exterior of a ball and hence the proofs in this paper can be applied in the case of black hole spacetimes. The assumption \cref{eq:SILED} is known to hold for black hole spacetimes.
\end{remark}

\begin{remark}
The case $n=2$ of \cref{eq:eqn}, namely the nonlinearity $c\phi^{2}\pa\p$, exhibits slow initial decay, so we prove an $r$-weighted integrated local energy decay estimate (see \cref{sec:rp}) to jumpstart the pointwise decay iteration.%
\end{remark}

\subsection{Outline of the paper} \label{ss:outline}
Here we overview the proof of \cref{thm:main} by presenting an outline of the paper. %
\begin{itemize}
\item
In \cref{sec:notation} we define notation that is used throughout the article.
\item
In \cref{sec:commuting} we make the transition from the uniform energy bounds and the integrated local energy decay statements of \cref{EB,eq:SILED}, which were stated only for derivatives, to their versions for vector fields. 
\item In \cref{sec:fromLEDtoptw}
we connect pointwise bounds to $L^2$ estimates and norms.
\item
In \cref{sec:preliminaries} we rewrite the equation in a way amenable to our pointwise decay iteration scheme. We state and prove lemmas that are used in the scheme to improve the pointwise decay rates of the solution. %
\item
In \cref{sec:rp} we prove an $r$-weighted local energy decay estimate which allows us to start the pointwise decay process for the slowly decaying nonlinearity $c\phi^2\pa\p$, i.e. the special case of \cref{eq:eqn} with the lowest possible power $n=2$. This section is unnecessary for all the other nonlinearities considered in this article. 
\item
In \cref{sec:ext} we prove the final decay rate for $\phi$ and its vector fields in the region exterior to the light cone $\{ r=t\}$, that is, in the region $\{ r \geq t \}$. 
\item
In \cref{sec:int} we prove the final decay rate for $\phi$ and its vector fields in the region inside of the light cone, that is, $\{ r \leq t\}$. 
\item
In \cref{sec:final} we demonstrate $\la t-|x|\ra\inv$ better decay for the solution to equations of the form \cref{eq:eqn,eq:eqn2}. 
\end{itemize}

\section{Notation}\label{sec:notation}

\subsection{Notation for dyadic numbers and conical subregions}
We work only with dyadic numbers that are at least 1. We denote dyadic numbers by capital letters for that variable; for instance, dyadic numbers that form the ranges for radial (resp. temporal and distance from the cone $\{|x|=t\}$) variables will be denoted by $R$ (resp. $T$ and $U$); thus $$R,T, U\ge 1.$$ 
	We choose dyadic integers for $T$ and a power $a$ for $R,U$---thus $R = a^k$ for $k\ge1$--- different from 2 but not much larger than 2, for instance in the interval $(2,5]$, such that for every $j\in\N$, there exists $j'\in\N$ with 
$a^{j'} = \f38 2^j.$
 
\subsubsection{Dyadic decomposition} Let $\R^{+} := [0,\infty)$. 
We decompose the region $\{r\le t\}$ based on either distance from the cone $\{r=t\}$ or distance from the origin $\{r=0\}$. We fix a dyadic number $T$ and we define the following dyadic sets
\begin{align*}
C^R_T &:=\begin{cases}
C_T\cap \{R<r<2R\} & R>1\\
C_T\cap \{0 < r < 2\} & R=1
\end{cases}\\
C^U_T &:=\begin{cases}
\{ (t,x) \in \R^{+} \times \rt  :  T\le t\le 2T\} \cap \{U<|t-r|<2U\} & U>1\\
\{ (t,x) \in \R^{+} \times \rt  :  T\le t\le 2T\} \cap \{0< |t-r|<2\} & U=1
\end{cases} %
\end{align*}
	We define
\begin{equation*}
C_T^{<3T/4} := \bigcup_{R < 3T/8} C_T^R.
\end{equation*}

Now letting $R > T$, we define
\begin{align*}
C^T_R &:= \{ (t,x) \in \R^{+} \times \rt  :  r \ge t, T \le t\le 2T, R \le r\le 2R, R \le |r-t| \le 2R\}
\end{align*} 

	$C_T^R, C_T^U,$ and $C^T_R$ are where we shall apply Sobolev embedding, which allows us to obtain pointwise bounds from $L^2$ bounds. 
Given any subset of these conical regions, a tilde atop the symbol $C$ will denote a slight enlargement of that subset on their respective scales; for example, $\crtt$ denotes a slightly larger set containing $\crt$.

\subsection{Notation for the symbols $n$ and $N$} \label{subsec:N}
Throughout the paper the integer $N$ will denote a fixed and sufficiently large positive number, signifying the highest total number of vector fields that will ever be applied to the solution $\p$ to \eqref{eq:problem} in the paper. 

We use the convention that the value of $n$ may vary by line.%

\subsection{The use of the tilde symbol} If $\xS$ is a set, we shall use $\ti\xS$ to indicate a slight enlargement of $\xS$, and we only perform a finite number of slight enlargements in this paper to dyadic subregions. The symbol $\ti\xS$ may vary by line. 

\subsection{Summation of norms}\label{conv:normsum}

Recall the subscript notation \cref{vf defn} for vector fields. Let $\| \cdot\|$ be any norm used in this paper. Given any nonnegative integer $N\ge0$, we write $\|g_{\le N}\|$ to denote $\sum_{|J|\le N} \|g_J\|$. For instance, taking the absolute value as an example of the norm, the notation $|\pm(t,x)|$ means
$$|\pm(t,x)| = \sum_{J : |J| \leq m} |\p_J(t,x)|.$$

\subsection{Other notation}
If $x =(x^1,x^2,x^3)\in\R^3$, we write 
\begin{align*}
u:= t-r,  \quad v := t+r.
\end{align*}

Next, we define the backward light cone in the (upper right quadrant of the) plane with apex $(r,t)$, which we denote as $\dtr$, and other related objects. Let 
$$\calR_1:= \{ R : R < u/8 \}, \quad\calR_2 := \{ R : u/8 < R < v \}, \qquad u > 0.$$
	Let $D_{tr}$ denote   
	\[
	D_{tr} := \{ (\rho,s)  \in \R_+^2: -(t+r) \leq s-\rho\leq t-r, \ |t-r| \leq s+\rho \leq t+r\}.
	\] 
When we work with $D_{tr}$ we shall use $(\rho,s)$ as variables, and $D_{tr}^{ R}$ is short for $D_{tr}^{\rho\sim R}$.
Thus
\begin{equation*}
\co:= D_{tr} \cap \{ (\rho,s) : R < \rho<2R\}, \ R>1; \qquad D^{R=1}_{tr} := \dtr\cap\{(\rho,s) : \rho<2 \}. 
\end{equation*}
We use the notation $dA := dsd\rho$ for integrations over $\dtr$ and its subsets.

\section{Vector field commutation}\label{sec:commuting}
We have
\begin{equation}\label{D decomp}
\pa w \in S^Z(\jr\inv) (\xO w, Sw) + S^Z(1) \pat w  \text{ if }r \ge t/2.
\end{equation}
This is clear for $\pa_t$ and the angular derivatives, while for $\pa_r$ we write
$\pa_r = r\inv (S- t \pat).$
	
	We define $\dot \calC$ to be the collection of real linear combinations of the operators
\begin{equation}\label{commdef}
 \pa s_{1+q'} \pa, \  s_{1+q'} \pa \pa, \ s_{2+q'} , \ \pa s_{1+q'},  \ s_{1+q'} \pa
\end{equation} 
where $q'>0$ is a number which depends on the assumptions made about the coefficients $h,g^\xo,V,A,$ and $B$ in \cref{thm:main}. Recall the counting convention \cref{cou.conv}. 
	We have the following fact:
\begin{lemma}
There are some operators $\dot C\in \dot\calC$ such that
\begin{equation}\label{rot and scal comm}
\xO^J (S+2)^k Pw = P \xO^J S^k w + \dot C w_{\le 4(|J|-1) + 10k}
\end{equation} where we interpret $\dot C w_{\le 4(|J|-1) + 10k}$ as a sum, and subscripts with negative real value denote the zero multiindex.

\end{lemma}

\noindent\textit{Sketch of proof.}
By the assumptions in the main theorem,
\begin{equation}\label{commfact1}
[P,\pa]  \in \dot\calC.
\end{equation}
\begin{equation}\label{commfact2}
[P,\xO]  \in \dot\calC.
\end{equation}
\begin{equation}\label{commfact3}
[P,S] - 2P - s_{2+\xs} \xO^2  \in \dot\calC.
\end{equation}
	We use \cref{commfact1,commfact2,commfact3} and the result follows by a proof by induction. Note that starting from $\xO^J(S+2)^k P$ and then commuting the vector fields with $P$, then other than $P\xO^JS^k$, the terms with the highest vector field count (assuming $g^\xo$ is not the zero function) are those of the form  
	$$\dot C (\xO,S)^{= |J| + k-1} w, \ \ \dot C\in \dot\calC;$$
	 more specifically, those of the form $\dot C \xO^{|J|-1}S^k$. This explains the subscript $4(|J|-1)+10k$ in \cref{rot and scal comm}. \hfill$\Box$

\begin{lemma}
Given the assumptions on the operator $P$, there exists a positive real number $q'>0$ such that for any multiindex $J$,
$$|P\p_J| \ls \f{|\p_{\leq |J|-1}|}{\jr^{2+q'}} + \f{|\pa \p_{\leq |J|}|}{\jr^{1+q'}} + |(P\p)_{\le |J|}|.$$
\end{lemma}

\begin{proof}
There is a constant $q'>0$ such that the operator $P$ can be written schematically as
$P=\Box+\pa s_{1+q'} \pa + s_{1+q'} \pao^2 + s_{2+q'} + s_{1+q'}\pa + \pa s_{1+q'}$. We have $[Z,\pa] = c\pa$ schematically, for some real number $c$ depending on $Z$. 	
	We include the terms arising from $g^\xo\Delta_\xo$ together with the $\jr^{-1-}|\pa \p_{\le|J|}|$ term.
\end{proof}

\begin{proposition}\label{led.for.vf}
Assume %
$$ \| \pa^{\le m}\p\|_{LE^1([T_0,T_1)\times\rt)} 
\ls_m  \|\pa \p(T_0)\|_{H^{m+k_0}(\rt)} + \|\pa^{\le m}(P\p)\|_{LE^*([T_0,T_1)\times\rt)}+\|\pat\pa^{\le m}\p\|_{LE[T_0,T_1)}.$$
Then we have
\begin{align}\label{prop.claim}
\begin{split}
\|\phi_{\le m}\|_{LE^1([T_0,T_1)\times\rt)} \ls \|\pa\p_{\le m+k_0}(T_0)\|_{\lt} + \|(P\p)_{\le m}\|_{(LE^*)([T_0,T_1)\times\rt)} + \|\pat\pm\|_{LE[T_0,T_1)}.
\end{split}
\end{align}
\end{proposition}

\begin{proof} We prove \eqref{prop.claim} by induction. The base case holds by the base case of \cref{eq:SILED}. Then
\begin{align*}
\|\p_{(I,J,k)}\|_{LE^1} 
    &\ls \|\pa \xO^JS^k\p(T_0)\|_{H^{|I|+k_{0}}} + \|\xO^JS^k(P\p)\|_{LE^{*,|I|}} + \|[P,\xO^JS^k]\p\|_{LE^{*,|I|}} \\
    &\qquad + \|\pat\pa^{\le|I|}\xO^JS^k\p\|_{LE}	\\
&\ls \|\pa\p_{\le m+k_{0}}(T_0)\|_{\lt} + \|(P\p)_{\le m}\|_{LE^*} + \|[P,\xO^JS^k]\p\|_{LE^{*,|I|}} + \|\pat\pm\|_{LE}\\
&\ls \|\pa\p_{\le m+k_{0}}(T_0)\|_{\lt} + \|(P\p)_{\le m}\|_{LE^*} + \|\jr^{-1-}\pa \p_{\le m-2}\|_{LE^*} + \|\jr^{-2-} \p_{\le m-2}\|_{LE^*} \\
&\qquad +\|\pat\pm\|_{LE}\\
&\ls \|\pa\p_{\le m+k_{0}}(T_0)\|_{\lt} + \|(P\p)_{\le m}\|_{LE^*} + \| \p_{\le m-2}\|_{LE^1}+\|\pat\pm\|_{LE}\\
&\ls \|\pa\p_{\le m+k_{0}}(T_0)\|_{\lt} + \|(P\p)_{\le m}\|_{LE^*}+\|\pat\pm\|_{LE}
\end{align*} %
 The final line follows by the induction hypothesis.
\end{proof}

\section{From local energy decay to pointwise bounds}\label{sec:fromLEDtoptw}

In this section we will show that local energy decay bounds imply certain slow decay rates for the solution, its vector fields, and its derivatives---see \cref{inptdcExt,derbound}. 

The following pointwise estimate for the second derivative will be used, for instance, when applying \cref{DyadLclsd} to the functions $ w = \pa\pm$ (that is, when we bound the first-order derivatives pointwise); this will be done in \cref{derbound}. 
\begin{lemma}\label{2ndDeBd'}
Assume $\p$ is sufficiently regular. Let $J_{n} := |J|+n$. Then for any point $(t,x)$
 \begin{equation}\label{2ndD}
|\pa^2\p_J(t,x)| \ls
\left(\f1\jr+\f1\ju\right)|\pa\p_{\le J_{n}}| +  \left( 1 + \f{t}\ju \right) \jr^{-2}|\p_{\le J_{n}}|
+ \left( 1 + \f{t}\ju \right) |(P\p)_{\le J_{n}}|.
 \end{equation}
\end{lemma}

\begin{proof}
First we note that 
\begin{equation}\label{r small}
\pa^2\p_J \ls \jr\inv |\pa\p_{\le |J|+n}|, \quad r \ls 1
\end{equation}
holds for any number $n \geq 1$. Thus \cref{2ndD} holds for all $r\ls 1$. 

Assume henceforth that $r \gg 1$. We have \begin{equation}\label{2nd'sEst}
|\pa^2\p_J| \ls
\left(\f1\jr+\f1\ju\right)|\pa\p_{\le|J|+n}|
+ \left( 1 + \f{t}\ju \right) |(\Box\p)_{\le|J|}|.
 \end{equation}
where the $|J|=0$ case follows from Lemma 2.3 in \cite{KS}. The general case follows after commuting with vector fields. By \eqref{P def} and \cref{coeff.assu} we have
\begin{align*}
(\Box - P)\p \in S^Z(\jr^{-1-\xs}) (\pa^2\p +\pa\p)
+ S^Z(\jr^{-2-\xs})\xO^{\le 2}\p  
\end{align*} 
Thus for some $m$ we have 
\begin{align}\label{error}
\begin{split}
\left|\left((\Box - P )\p\right)_{\leq |J|}\right| 
&\in S^Z(\jr^{-1-\xs})|\pa^2\p_{\leq |J|+m}| + S^Z(\jr^{-2-\xs}) |\p_{\leq |J|+m}| \\
&\in S^Z(\jr^{-1-\xs}) |\pa\p_{\leq |J|+n}| + S^Z(\jr^{-2-\xs}) |\p_{\leq |J|+n}|, \quad n = m+1.
\end{split}
\end{align}
	The conclusion now follows from \cref{r small,2nd'sEst,error}. 
\end{proof}

The primary estimates that let us pass from local energy decay to pointwise bounds are contained in the following lemma. 
\begin{lemma}
\label{DyadLclsd}
Let $w\in C^4$,  $Z_{ij} :=S^i\xO^j$,  $\mu:= \la \min(r,|t-r|) \ra$, and $\calR \in \{\crt,\cut,\ctr\}$. 
	Then we have  \begin{equation}\label{DyadLclsdBd}
 \| w\|_{L^\iy(\calR)} \ls \sum_{i\le 1,j\le 2} \f1{|\calR|^{1/2}}  \lr{  \|Z_{ij} w\|_{L^2(\calR)} +  \|\mu \pa Z_{ij} w\|_{L^2(\calR)} }.
\end{equation}
where we assume $1\ll U \le \f38 T$, $1\ll R \le \f38 T$ and $R > T \gg1$ in the cases $\cut,\crt,\ctr$ respectively, and $|\calR|$ denotes the measure of $\calR$. 
\end{lemma}

\noindent \textit{Sketch of proof.} One uses exponential coordinates, which results in $\calR$ being transformed into a region of constant size in all directions. Then one uses the fundamental theorem of calculus for the $s,\rho$ variables and Sobolev embedding for the angular variables. Finally, changing coordinates to return to the original region $\calR$ produces the $|\calR|^{-1/2}$ factor. \hfill $\Box$

The next proposition yields an initial global pointwise decay rate for $\phi_J$ under the assumption that the local energy decay norms are finite. We shall improve this rate of decay in future sections (see \cref{sec:ext,sec:int}) for solutions to \eqref{eq:problem}, culminating ultimately in the final pointwise decay rate stated in the main theorem. 

\begin{proposition} \label{inptdcExt}
Let $T$ be fixed and $\p$ be any sufficiently regular function.
There is a fixed positive integer $k$, %
such that for any multi-index $J$ with $|J|\le N - k$, we have:
\begin{align}\label{u/v decay}
\begin{split}
|\p_{J}| \leq \bar C_{|J|} \|\p_{\leq |J|+k}\|_{LE^1[T, 2T]} \ju^{1/2}\jv\inv. %
\end{split}\end{align} %
\end{proposition}
\noindent\textit{Sketch of proof.} We only sketch the proof here; full details are provided in Lemma 4.1 of \cite{L}. One uses \cref{DyadLclsd}, which proves \cref{u/v decay} except in the wave zone. For the wave zone, an extra Hardy-like inequality, which is proven by multiplying by a cutoff function localised to the wave zone, is used. It is written as Lemma 4.3 in \cite{L}. \hfill$\Box$

\subsection{Derivative bounds}
\begin{lemma}\label{lem:derbd}
If $w$ is a sufficiently smooth function, then for $\calR\in\{\crt,\cut,\ctr\}$ and for $R\gg1, U\gg1, T\gg1$ respectively,
\begin{equation} \label{easy}
  \| \pa w\|_{\lt(\calR)} 
  	\ls  \|\mu\inv w_{\le1}\|_{\lt(\ti\calR)} +  \|\jr P w\|_{L^2(\tilde \calR)}
\end{equation}

\end{lemma}

\begin{proof} We only prove the case $\crt$ as the other two cases are similar; see also Lemma 5.1 in \cite{L}.
Let $\chi(t,r)$ be a radial cutoff function on $\R^{1+3}$ with $\supp\chi \subset \crtt$ and $\chi=1$ on $\crt$. Note that
\begin{enumerate}
\item
If $r < t$ then for a sufficiently large constant $C'$, we have %
\begin{equation}\label{1f}
\chi \lr{ \f{u}t |\de w(t,x) |^2  }\le \chi\lr{   |\grad w |^2- w _t^2+ \f{C'}{ut} |S w |^2 }
\end{equation}
as an expansion of the terms $|Sw|^2, |\de w|^2$ reveals.%

\item Integrating by parts,
\begin{align}\label{ini.com}
\int \chi ( |\grad w |^2 -  w _t^2) \, dxdt
	&=  \int \chi w (\pat^2-\Delta) w \,dxdt - \int \f12 (\pat^2-\Delta)\chi  w ^2\, dxdt.
\end{align} There are no boundary terms in either time or space because of the compact support of $\chi(t,r)$ in both time and space. 
\end{enumerate}
 	Integrating \cref{1f} in spacetime, we have via \cref{ini.com}
\begin{align}\label{1g}
\int\chi \f{u}t |\de w |^2 \,dxdt
	&\le \int \chi w(\pat^2-\Delta) w + O(|\Box\chi| w ^2) + \f{C'}{ut}\chi |S w |^2 \,dxdt.
\end{align}
The proof of \cref{easy} will be complete once we incorporate $Pw$ into \cref{1g}:
\begin{itemize}\label{list}
\item
For $\int (\chi w)\paa(\hab \pab w) \,dxdt$, we integrate by parts and use Cauchy-Schwarz. A term
$$\int \chi \hab \paa\p \pab \phi \,dxdt = O\lr{ \int \chi \f{ |\de w|^2}\jr \,dxdt }$$
arises, and for this term we use the hypothesis that $R \gg 1$ for $h \neq 0$. 

Similarly, $$\int (\chi w)  (g^\xo\Delta_\xo w) \,dxdt$$ is treated by integration by parts and Cauchy-Schwarz. We use the smallness of $\jr^{-2-\xd'}=O(\jr\inv)$ for sufficiently large $R$.
\item
We use the bound $V\ls \jr^{-2}$.

\item
For $$\int \chi w B \pa w \,dxdt$$ we use Cauchy-Schwarz. The bounds we obtain are sufficient to prove the claim \cref{easy} when $B\in S^Z(\jr^{-1})$ in this part.
\end{itemize}

Assuming $\Box\chi\ls \jr^{-2}$, separating
$|\chi wPw| \ls \chi [ (R\inv w)^2 + (RPw)^2 ]$ in the right-hand side of 
\cref{1g}, and using the reasoning in the bullet points (along with the triangle inequality) to deal with $$\int (\chi w) ( (\Box-P) w) \,dxdt,$$
this proves the desired claim \cref{easy}.
\end{proof}

\begin{corollary}\label{cor:derbd}
If $w$ is a sufficiently smooth function, then for $\calR\in\{\crt,\cut,\ctr\}$
\begin{equation} \label{easy}
  \| \pa w_{\le m}\|_{\lt(\calR)} 
  	\ls  \|\mu\inv w_{\le m+n}\|_{\lt(\ti\calR)} +  \|\jr (P w)_{\le m}\|_{L^2(\tilde \calR)}
\end{equation}
\end{corollary}

\begin{proof}
Compared to the proof of \cref{lem:derbd}, here one only has to bound $$\int \chi w_{\le m} [P,Z^{\le m}] w \,dxdt.$$ 
Similar arguments involving integration by parts and Cauchy-Schwarz as those seen in \cref{lem:derbd} establish the claim.
\end{proof}

\begin{lemma}\label{lem:Sob}
Given a function $f : \R^{3}\to \R$, we have 
$$\|f\|_{L^{\iy}(R < |x| < R+1)} \ls R\inv \|f_{\le 3}\|_{\lt(R - 1 < |x| < R+2)}.$$

Given a radial function $f:\R^{3}\to\R$, we have
\begin{equation}\label{bound}
f\ls \jr^{-1/2} \|(f_{0},f_{1})\|_{\dot H^{1}\times \lt}
\end{equation}
where $f(0,x) = f_{0}(x),  \pat f(0,x) = f_{1}(x).$
\end{lemma}

\noindent \textit{Sketch of proof.} 
The first result can be proven by combining a localised embedding and an embedding on $S^{2}$.  \hfill $\Box$

Thus\begin{equation}\label{cor:Sob}
\pa \pm \ls \jr\inv C_{0}
\end{equation}
by the uniform boundedness of energy, where $C_{0}$ is a constant depending on the initial data. We also gain access to using \cref{bound} for $\p$ by way of the positivity of the fundamental solution of the wave operator in three spatial dimensions.

The next proposition shows that the first-order derivative (of solutions to \eqref{eq:problem}) decays pointwise faster by a rate of $\min(\jr,\nm)$. It utilises the initial global decay rate \cref{u/v decay}. 
\begin{proposition}\label{derbound}
Let $\p$ solve \eqref{eq:problem}, and assume that 
$$\pmn \ls \jr^{-\x}\jt^{-\xb}\ju^{-\eta}$$
for some sufficiently large $n$. We then have
\begin{equation}\label{claim}
\pa\pm \ls \jr^{-\x}\jt^{-\xb}\ju^{-\eta} \mu\inv, \quad \mu := \la \min(r,| t-r |) \ra.
\end{equation}
For solutions to the equation \cref{eq:eqn}, \cref{claim} also holds if $c = c(t,x)$ is non-constant if we assume that $\pa c \in S^{Z}(\jr^{-1})$.	
\end{proposition}

\begin{proof} Let $\calR \in \{ \cut, \crt, C^T_R\}$. We first prove the result for all the nonlinearities aside from $\phi^{2}\pa\phi$.

Given a function $w$, we have by \cref{lem:derbd}
\begin{equation}\label{1stDeBd}
\|\pa w_{\le m}\|_{\lt(\calR)} \ls \| \f{ w_{\le m+n}}{\mu}  \|_{\lt(\ti\calR)} + \|\jr (Pw)_{\le m} \|_{\lt(\ti \calR)}.
\end{equation} 
By \cref{u/v decay,cor:Sob} and the assumptions on all of the nonlinearities $P\p$ except for $\p^{2}\pa\p$,
\begin{equation}\label{1stDeBd'}
\|\pa\pm\|_{\lt(\calR)} \ls \| \f\pmn\mu \|_{\lt(\ti\calR)}.
\end{equation}
	By \cref{DyadLclsd}, we have 
\begin{align}\label{final.calc}
\begin{split}
\|\pa\pm\|_{L^\iy(\calR)} 
&\ls |\calR|\invh \sum \|Z\pa \pm\|_{\lt(\calR)} + \|\mu \pa Z \pa \pm \|_{\lt(\calR)} \\
&\ls |\calR|\invh \lr{ \|\pa\pmn\|_{\lt(\calR)} + \|\mu \pa^2 \pmn\|_{\lt(\calR)}	} \\
&\ls |\calR|\invh \lr{ \|\mu\inv \pmn\|_{\lt(\ti\calR)} + \|\mu \pa^2 \pmn\|_{\lt(\calR)}	} \\
&\ls |\calR|\invh \lr{ \|\mu\inv \pmn\|_{\lt(\ti\calR)} + \| \mu \left( \f1\mu|\pa\pmn| + ( 1 + \f{t}\ju ) \jr^{-2}|\pmn| \right) \|_{\lt(\calR)}	} \\
&\ls |\calR|\invh \lr{ \|\mu\inv \pmn\|_{\lt(\ti\calR)} + \|\mu(1 + \f{t}\ju)\jr^{-2} \pmn \|_{\lt(\calR)} }\\
&\ls |\calR|\invh \|\mu\inv \pmn\|_{\lt(\ti\calR)} 
\end{split}
\end{align}
which follows by \cref{1stDeBd'}. The final line follows because $\mu^2(1 + t/\ju) \ls \jr^2$. \ Finally, the claim \cref{claim} follows because $$\|\mu\inv \pmn\|_{\lt(\ti\calR)} \ls |\calR|^{1/2} \|\mu\inv\pmn\|_{L^\iy(\ti\calR)}.$$

In the case $P\p=c\p^{2}\pa\p$, we return to \cref{ini.com}. Instead of using Cauchy-Schwarz here we use the structure of the nonlinearity to integrate by parts to observe that
\begin{align}\label{the.}
\begin{split}
\int \chi \p \cdot \p^{2}\pa\p \,dxdt &=- \int\f13 \p^{3} \cdot \pa(\chi\p) \,dxdt\\
	&\ls \int \p^{6}\chi + |\pa\chi| \p^{4} \,dxdt + \eps' \int \chi (\pa\p)^{2}\, dxdt
\end{split}
\end{align}for some small $\eps'>0$, and then Cauchy-Schwarz and \cref{bound} (applied to the $\p^{6}$ term to obtain $\ls \jr\inv \p^{4}$) imply the desired claim \cref{1stDeBd'} and the rest of \cref{final.calc} follows. We may use \cref{bound} because of the positivity of the fundamental solution of $\Box$. %
For the non-constant $c(t,x)$ case, the hypothesis $\pa c\in S^{Z}(\jr\inv)$ combined with this calculation \cref{the.} shows the result as it is essentially identical to the $|\pa\chi|\p^{4}$ term with only slightly immaterially different support properties.  \ Taking \cref{1f} into account, slightly weaker bounds on the derivative $\pa c$ than $\pa c \in S^Z(\jr\inv)$\footnote{as if we use \cref{u/v decay} to bound two copies of $\p$, then
$\p^{4}|\pa c|\chi \ls \chi \ju\jv^{-2}|\pa c| \p^{2}$ and we then require $\f{t}u  \chi \ju\jv^{-2}|\pa c| \p^{2}\ls\chi \mu^{-2}\p^{2}$ where $\mu = \la \min(r,|u|)\ra$, so 
 $\pa c \in S^Z(\f{\jv^2}{\mu^{2}t})$ works. This bound allows $\jt$ growth along the light cone $\{r=t\}$ and also along the hypersurface $\{ r = 0 \}$, so it is much weaker than the $\pa c\in S^{Z}(\jr\inv)$ assumption there.
 } 
 are possible, however for simplicity we simply require this. 
\end{proof}

\section{Preliminaries for the iteration}\label{sec:preliminaries}
\begin{remark}[The initial data] \label{rem:id}
Let $w := S(t,0)\p[0]$ denote the solution to the free wave equation with initial data $\phi[0]$ at time 0. Then for any $|J| = O_N(1)$,
\begin{equation}\label{Kir}
w_J(t,x) = \f1{|\pa B(x,t)|} \int_{\pa B(x,t)} (\p_0)_J(y) + \nabla_y (\p_0)_J(y) \cdot (y-x) +t (\p_1)_J(y) \, dS(y).
\end{equation}
By \cref{Kir} and the assumptions $\p_0 \in \lt(\R^3)$, 
$$\|\jr^{1/2 + \xk} \pa\p_{\le N}(0)\|_{\lt(\rt)} <\iy, \quad \xk = \min( \sigma, \calT + \calN -3 )$$
we have
$$ w_J \ls \jv\inv\ju^{-1-\xk} .$$
\end{remark}

\subsection{Summary of the iteration} \label{outline of iteration}
By \cref{rem:id}, we may assume zero initial data in the following iteration. Second, note that it suffices to prove bounds in $|u| \geq 1$, because the desired final decay rate in $|u| < 1$ already holds by \cref{u/v decay}. Third, we distinguish the nonlinearity and the coefficients of $P - \Box$, and for both of these, we apply the fundamental solution. We iterate these two components in lockstep with one another. 

Due to the domain of dependence properties of the wave equation, we shall first complete the iteration in $\{ u < -1 \}$. For the iteration in $\{ u > 1\}$, the decay rates obtained from the fundamental solution are insufficient in the region $\{ r < t/2\}$, so we prove \cref{convrsn}. With the new decay rates obtained from \cref{convrsn}, we are then able to obtain new decay rates for the solution and its vector fields. At every step of the iteration, \cref{conversion} is used to turn the decay gained at previous steps into new decay rates. 

\begin{remark}[Reduced and irrational $\xs$] \label{smllirr.xs}
To simplify the iteration, we shall reduce the value of $\sigma$ if necessary to be equal to some positive irrational number less than the original value of $\xs$. %
We do this to avoid the appearance of logarithms in the iterations for $\phi_1$ and $\phi_2$ (see the decomposition \cref{decomp} below). 
We take $0< \xs \ll 1$. %
	\, In the sections spelling out the details of the iteration, namely \cref{sec:ext,sec:int}, we explain how we reach the final decay rate in \cref{thm:main} (wherein the \textit{original} value of $\xs$ is included in the final decay rate). 
\end{remark}

\subsection{Setting up the problem}\label{settingup}

We rewrite \eqref{eq:problem} as
\[
\Box\p = (\Box - P)\p + F = -\pa_\alpha(h^{\alpha\beta}\pa_\beta\p +B^{\alpha}\p) - g^\xo \Delta_\xo \p - (V-\pa_\alpha B^{\alpha})\p + F, \quad F := P\p
\]

Using the assumptions \eqref{coeff.assu}, we can write this as
$$\Box\p \in \pa \left(S^Z(\jr^{-1-\xs}) \p_{\leq 1}\right) + S^Z(\jr^{-2-\xs}) \p_{\leq 2} + F$$
	After commuting with the vector field $Z^{J}$,
\begin{equation}\label{first write}
\Box\p_J \in \pa \left(S^Z(\jr^{-1-\sigma}) \p_{\leq m+1}\right) + S^Z(\jr^{-2-\sigma}) \p_{\leq m+2} + F_{\le m}, \quad m=|J|
\end{equation}
	Due to the derivative gaining only $\ju\inv$ in the wave zone, as was proved in \cref{derbound}, we shall perform a decomposition as follows.
First, we note that, for any function $w$, 
\begin{equation}\label{D decomp}
\pa w \in S^Z(\jr^{-1}) w_{\leq 1} + S^Z(1) \pat w, \quad r\geq t/2
\end{equation}
which is clear for $\pa_t$ and $\pa_\xo$, while for $\pa_r$ we write
$\pa_r = \frac{S}{r}- \frac{t}{r}\pa_t.$
	Let $\chi_\text{cone}$ be a cutoff adapted to the region $\{t/2 \le r \le 3t/2\}$. We now rewrite \eqref{first write} as
\begin{equation}\label{final write}
\Box\p_J \in S^Z(\jr^{-2-\sigma}) \p_{\leq m+2} + (1- \chi_\text{cone}) \left(S^Z(\jr^{-1-\xs}) \pa\p_{\leq m+1}\right) + \pa_t \left(\chi_\text{cone} S^Z(\jr^{-1-\sigma}) \p_{\leq m+1}\right) + F_{\le m}
\end{equation}
	We now write $\p_J = \sum_{j=1}^3\p_j$ where 
\begin{equation}\label{decomp}
\begin{split}
\Box \p_1 = G_1, \quad G_1 \in S^Z(\jr^{-2-\sigma}) \p_{\leq m+2} + (1- \chi_\text{cone}) \left(S^Z(\jr^{-1-\sigma}) \pa\p_{\leq m+1}\right) \\
\Box \p_2 = \pa_t G_2, \quad G_2\in \chi_\text{cone} S^Z(\jr^{-1-\sigma}) \p_{\leq m+1} \\
\Box \p_3 = F_{\le m} = G_3.
\end{split}
\end{equation}

In the sections for the pointwise decay iteration (\cref{sec:int,sec:ext}) we will be comparing our nonlinearity $P\p$ in \cref{eq:problem} to the case of the semilinear wave equation satisfying the null condition. Thus we record the following observation:
\begin{align}\label{tang.prop}
\begin{split}
(\pat+\pa_{r}) \p &= \frac{\mu u}{v'} \pa \p + \frac{1}{v'} S\p, \qquad v' \in \{ r,t\}, \ \mu \in \{-1,1\} \\
\pao \p &\in S^{Z}(\la r\ra\inv) \xO\p.
\end{split}
 \end{align}
 We will only be using the following simpler form:
 \begin{align*}
\begin{split}
(\pat+\pa_{r}) \p &= \frac{\mu u}{v'} \pa \p + \frac{1}{v'} Z\p, \qquad v' \in \{ r,t\}, \ \mu \in \{-1,1\} \\
\pao \p &\in S^{Z}(\la r\ra\inv) Z\p.
\end{split}
 \end{align*}

\subsection{Estimates for the fundamental solution}\label{sec:estsfdmt}
\begin{lemma}\label{conversion}
Let $m\ge0$ be an integer and suppose that $\psi : [0,\iy)\times\R^3\to\R$ solves $$\Box\psi (t,x)= g(t,x), \qquad \psi(0) = 0, \quad \pa_t \psi(0) = 0. $$ 
Define
\begin{equation}\label{hdef}
h(t,r) := \sum_{i=0}^2 \|\Omega^i g (t, r\omega)\|_{L^2(\mathbb{S}^2)}
\end{equation}
	Assume that %
$$ h(t,r)  \ls \f{1}{ \jr^\x \la v\ra^{\beta} \la u\ra^\eta }, \quad \x \in (2,3) \cup (3,\iy) , \quad \beta\geq 0, \quad \eta\geq -1/2.$$
	Define%
\[
\tilde\eta = \left\{ \begin{array}{cc} \eta -2,& \eta<1   \cr -1, & \eta > 1 
  \end{array} \right. .
\]
	We then have in both $\{ u > 1\}$ (without any additional restrictions on the value of $\x+\xb+\eta$), and $\{u < -1\}$ in the case $\x+\xb+\eta>3$:
\begin{equation}\label{Bd1}
\psi(t, x)\lesssim \frac{1}{\la r\ra\la u\ra^{\alpha+\beta+\tilde\eta-1}}.
\end{equation}

On the other hand, if $\alpha+\beta+\eta < 3$ and $u < -1$, we have
\begin{equation}\label{Bd2}
\psi(t, x)\lesssim r^{2-(\alpha+\beta+\eta)}. 
\end{equation}

\end{lemma}

\begin{proof}
A detailed proof of \eqref{Bd1} can be found in Lemma 6.5 of \cite{L}. The idea of the proof is to use Sobolev embedding and the positivity of the fundamental solution of $\Box$ to show that
$$r\psi \ls \int_{D_{tr}} \rho h(s,\rho) ds d\rho,$$
where $D_{tr}$ is the backwards light cone with vertex $(r, t)$, and to use \cref{hdef}.

Let us now prove \eqref{Bd2}, which was subject to the hypotheses $\x+\xb+\eta<3$ and $u<-1$. In this case ${D_{tr}} \subset \{r-t\leq u'\leq r+t, \quad r-t\leq \rho\leq r+t\}$ and we obtain, using that $ u' \le \rho$ and $\rho \gs \rho+s$ in $D_{tr}$:
\[
r\psi \lesssim \int_{r-t}^{r+t} \int_{u'}^{r+t} \jrho^{1-\alpha-\beta} d\rho \, \la u'\ra^{-\eta} du' \lesssim \int_{r-t}^{r+t} \la u'\ra^{2-(\alpha+\beta+\eta)} du' \lesssim (t+r)^{3-(\alpha+\beta+\eta)}
\] where the final bound follows from the hypothesis that $\x+\xb+\eta<3$. 
This finishes the proof because $t+r \le 2r$ when $u<-1$. 
\end{proof}

For the function $\phi_2$, which was written down in order to deal with the metric terms near the light cone, we will use the following result for an inhomogeneity of the form $\pat g$ supported near the cone. The result is similar to \cref{conversion}, aside from a gain of $\ju$ in the estimate: see \cref{Bd1der}.

\begin{lemma}\label{Minkdcyt}
Let $\psi$ solve 
\begin{equation}\label{Mink2}
\Box \psi = \pa_t g, \qquad \psi(0) = 0, \quad \pa_t \psi(0) = 0,
\end{equation}
where $g$ is supported in $\{ 1/2 \leq |x|/t \leq 3/2 \}$. Let $h$ be as in \cref{hdef}, and assume that 
\[
|h| + |Sh| + |\Omega h| + \la t-r\ra |\pa h| \ls \frac{1}{\la r\ra^{\alpha}\la u\ra^{\eta}},  \quad 2 < \alpha <3,  \quad \eta\geq -1/2.
\]
	Then in $\{ u > 1\}$, and $\{ u < -1\}$ when $\alpha+\eta > 3$
\begin{equation}\label{Bd1der}
\psi(t, x)\lesssim \frac{1}{\la r\ra\la u\ra^{\alpha+\tilde\eta}}.
\end{equation}

\end{lemma}
\begin{proof}
Let $\tpsi$ solve $$\Box \tpsi = g, \quad \tpsi(0)= 0, \quad \pat \tpsi(0)=0.$$
In the support of $g$ we have
\[
(t \partial_i + x_i \partial_t) h \ls |Sh| + |\xO h| + \la t-r\ra |\pa_r h|.
\]	
	By \cref{conversion} with $\beta=0$ applied to $\pa \tpsi$,  
$\Omega \tpsi$, $S \tpsi$, and the bound
$$\la u\ra \partial_t \tpsi \ls |\pa \tpsi| +|S\tpsi|+|\xO \tpsi|  + \sum_i  | (t \partial_i + x_i \partial_t) \tpsi|$$
the claim follows.
\end{proof}

\section{Preliminaries for the nonlinearity $\p^{2}\pa\p$}\label{sec:rp}
In the following theorem, we prove an $r$-weighted integrated local energy decay statement for solutions of the equation $P\p = \sum_{i=1}^{M}F_{i}(\p,\pa\p)$ where each $F_{i}(\p,\pa\p)$ is of the form $c(t,x)\p^{2}\pa\p$ for some function $c(t,x) \in S^{Z}(1)$.

\begin{theorem}[The $r^\xg$ estimate] \label{the rp est} Let $\p$ %
$$P\p=\sum_{i=1}^{M} c_{i}(t,x) \phi^{2}\pa_{(i)}\p, \quad \pa_{(i)}\in \{ \paa \}_{\x \in \{0,\dots,3\}}, \quad c_{i}\in S^{Z}(1).$$
Let $\gamma < 2\sigma, \gamma < 1$.
 Let $T_2 > T_1 \geq 0$.
Fix $m \in \N_{0}$. Assume that for a sufficiently large $n$, we have
\begin{equation}
\label{small leo}
\|\pmn\|_{LE^{1}}\ll1.
\end{equation}
 Then
\begin{equation}\label{rg.est}
A_{\xg,m}+ E^\xg_{\pm}(T_2)\ls E^\xg_{\pm}(T_1) +  \|\pa \p_{\le m}\|_{LE(T_1,T_2)}^2 + \|\pa^2\phi_{\le m}\|_{LE(T_1,T_2)}^2
\end{equation}
where the $A,E$ norms are:
$$A_{\xg, m} := \int_{T_1}^{T_2}\int_{\R^3} (\pm)^2 r^{\xg-3} + |\bar\partial\pm|^2 r^{\xg - 1} \, dxdt$$
$$E^\xg_{\pm}(T_1) := \| r^{\xg/2} ( \pao\pm,(\pav+\f1{2r})\pm,\f\pm{r} ) (T_1)\|^2_{L^2(\R^3)}, \  \| r^\x (f_1,\dots,f_n)\| := \sum_{j=1}^n\|r^\x f_j \|.$$
\end{theorem}

\begin{proof} 
Fix $m\ge0$. Let $|J| \leq m$. Fix $0 \leq T_1 < T_2$. Let
$$A_{\xg, J} := \int_{T_1}^{T_2}\int_{\R^3} \p_J^2 r^{\xg-3} + |\bar\partial\phi_J|^2 r^{\xg - 1} dxdt.$$
 
Integrating by parts in the region $[T_{1},T_{2}] \times \R^{3}$,
\begin{align}\label{Jmain}
	\begin{split}
	\iint &\Box \p_J(r^\gamma\pav \p_J + r^{\gamma-1} \p_J) \,dx\,dt = 
	\iint  - \f{\gamma r^{\gamma-1}}2(\pav\p_J)^2 - \f12(2 - \gamma)r^{\gamma-1}|\pao\p_J|^2 \\
	&-\f{\xg(1-\xg)r^{\xg-3}}2\p_J^2 \, dxdt +\int_{\R^3}- r^\xg\left[ \f12|\pa \p_J|^2 + \pa_r\p_J\pat\p_J + \f1{2}\f{\p_J}{r}\pat\p_J\right]_0^T \,dx
	\end{split}	
	\end{align}
\begin{itemize}
\item
We have
\begin{align*}
\int_{\R^3} -r^\gamma \f1{2}\f{\p_J}{r}\pat{\p_J} \,dx 
	&= \int -r^\gamma \f1{2}\f{\p_J}{r}(\pav-\pa_r){\p_J} \,dx  \\
	&= \int -r^\gamma \f1{2}\f{\phi_J}{r}\pav\p_J \, dx + \int_{S^2}\int_0^\iy  r^\xg \f12 \f{\p_J}{r} \pa_r{\p_J} \, r^2 \, drd\xo  \\
	&= \int -r^\gamma \f1{2}\f{\p_J}{r}\pav\p_J\, dx -\int_{S^2}\int_0^\iy \f{\gamma+1}4 r^\gamma{\p_J}^2 \, drd\xo \\
	&= \int -r^\gamma \f1{2}\f{\p_J}{r}\pav\p_J\, dx -\f{\xg+1}4 \int_\rt r^{\gamma}\f{{\p_J}^2}{r^2} \,dx
\end{align*}

Thus,
\begin{align}\label{bdry}
\begin{split}
- \int_\rt &r^\gamma\left( \f12|\pao\phi_J|^2  + \f12 (\pav\phi_J)^2 + \f{\gamma+1}4 \f{\phi_J^2}{r^2} + \f12 \f{\phi_J}{r}\pav\phi_J\right)_{T_1}^{T_2} \,dx \\
	&=- \int_\rt r^\gamma \left( \f12|\pao\phi_J|^2 + \left[\f12(\pav\phi_J)^2 + \f18\f{\phi_J^2}{r^2} + \f12 \f{\phi_J}{r}\pav\phi_J \right] + ( \f\gamma4 + \f18 ) \f{\phi_J^2}{r^2}  \right)_{T_1}^{T_2} \,dx \\
	&=-\int_\rt r^\gamma \left( \f12|\pao\phi_J|^2 + \f12\left[ \pav\phi_J + \f{\phi_J}{2r} \right]^2 + \left( \f\gamma4 + \f18 \right) \f{\phi_J^2}{r^2}  \right)_{T_1}^{T_2} \,dx
\end{split}
\end{align}

\item In this item, we deal with the nonlinearity $P\p = \phi^{2}\partial\phi$. By \cref{u/v decay}, the assumption of small $LE^{1}$ norm in \cref{small leo}, \cref{derbound}, as well as Cauchy-Schwarz, we obtain
\begin{align}
\begin{split}
\int P\p (r^{\xg}\pav\p+r^{\xg-1}\p) \,dxdt &\ls | \int \p^{2}\pa\p \cdot r^{\xg}(\pav\p + r\inv\p) \,dxdt | \\
&\ls \eps' A_{\xg,m} 
\end{split} 
\end{align}for a sufficiently small number $\eps'>0$ depending on the initial data, which allows us to absorb $\eps' A_{\xg,m}$ to the left-hand side.  Here the bound \cref{bound} would not have sufficed, and instead we used \cref{u/v decay}.

    \item 
   \begin{enumerate}
    \item Let  $V\in S^Z(\jr^{-2-})$. We split the following integral into small-$r$ and large-$r$ regions. In the large-$r$ region we use the positivity of $\xs$ to obtain a small coefficient. 
\begin{align}\label{potential}
\begin{split}
\int_{T_1}^{T_2} \int_{\R^3} &|V_{\le m} \pm r^\xg (\pav \p_J + \f{\p_J}{r}) | \,dxdt 
   \ls  \int r^{\xg} \f1{\jr^{2+\xs}}|\pm| ( |\pav \p_J|+|\f{\p_J}r| ) \,dxdt \\
   &\le \int_{r\gg1} \f1{\jr^{\xs}} r^{\xg-2}|\pm| ( |\pav \p_J|+|\f{\p_J}r| ) \,dxdt + \int_{r\ls1} r^{\xg} \f1{\jr^{2+\xs}}|\pm| ( |\pav \p_J|+|\f{\p_J}r| ) \,dxdt \\
   &\ls \eps_{\xs} \int_{r \gg1} \f{(\pm)^2 + \p_J^2}{r^{3-\xg}}+ \f{|\pav\p_J|^2}{r^{1-\xg}}  \, dxdt + \int_{r\ls1} r^{\xg} \f1{\jr^{2+\xs}}|\pm| ( |\pav \p_J|+|\f{\p_J}r| ) \,dxdt 
   \\
   &\ls \eps_{\xs} A_{\xg, m} + \|\pm\|_{\leo(T_{1},T_{2})}^{2}
\end{split}
\end{align}
where $\eps_{\xs}>0$ is a sufficiently small constant. We bounded the small-$r$ integral by the usual ILED norms $\|\pm\|_{\leo}.$

If $B \in S^Z(\jr^{-1-\xs_B})$ and $2 \sigma_B > \gamma$:
\begin{align}\label{B term}
\begin{split}
\int_{T_1}^{T_2} \int_{\R^3} |B_{\le m} \pa\pm r^{\xg-1}\p_J| dxdt
  &\ls  \int \f1{\jr^{1+\xs_B}} |\pa\pm r^{\xg-1}\p_J|  \\
  &\ls \f1\eps\int \f{|\pa\pm|^2}{\jr^{1 + 2\xs_B - \xg}} + \eps \int r^{\xg-3}\p_J^2 \\
  &\ls \f1\eps \|\pm\|_{LE^1(T_1,T_2)}^2 +  \eps A_{\xg, J}
\end{split}
\end{align}  where $\eps$ is a small constant.
	The bound on $\int |B_{\le m}\pa\pm r^\xg| \cdot |\bar\pa\phi_J|dxdt$ is similar. 

\item We may schematically write all terms involving the metric $\hab$ as $\int (|\pa h_{\le m} \pa\pm| + |h_{\le m} \pa^2\pm|) r^\xg(\f{|\p_J|}r+|\pav\p_J|)dxdt$, where $|J| = m$. Here we assume only that $h\in S^{Z}(\jr^{-1-\xs})$. 
\begin{align}\label{h term}
\begin{split}
 \int_{T_1}^{T_2} \int_{\R^3} &(|\pa h_{\le m} \pa\pm| + |h_{\le m} \pa^2\pm|) r^\xg(\f{|\p_J|}r+|\pav\p_J|)dxdt \\
	&\ls  \int \f1{\jr^{1+\xs}}r^\xg\left(|\pa\pm| + |\pa^2\pm|\right) (|r\inv \p_J| + |\pav\p_J|) \,dxdt\\
	&\ls  \int \f{r^{\f{\xg+1}2}}{\jr^{1+\xs}} \left(|\pa\pm| + |\pa^2\pm|\right)  r^{\f{\xg-1}2} (|r\inv \p_J| + |\pav\p_J|) \,dxdt    \\
	&\ls \f1\eps \int \f1{\jr^{2+2\xs}}r^{\xg+1}  \left(|\pa\pm| + |\pa^2\pm|\right)^2 \,dxdt +\eps A_{\xg, J} \\
	&\ls  \|\pa\pm\|_{LE(T_1,T_2)}^2 +  \|\pa^2\pm\|_{LE(T_1,T_2)}^2 + \eps A_{\xg, J}	    \text{ if } 2 \xs > \xg
\end{split}
\end{align} where $\eps$ is a small constant. We absorb $\eps A_{\xg,J}$ to the left-hand side.
    
\end{enumerate}
    
\end{itemize}

Taking the sum of \cref{Jmain,bdry,potential,B term,h term} over all $|J| \leq m$, i.e. $\sum_{|J| \leq m}$ (\cref{Jmain,bdry,potential,B term,h term}), 
  we get
$$A_{\xg, m} + E^\xg_{\pm}(T) \ls E^\xg_{\pm}(0) +  \|\pm\|^2_{LE^1(T_1,T_2)} +  \|\pa^2 \pm\|_{LE(T_1,T_2)}^2.$$
\end{proof}

Assume that $u >0$. The next lemma shows that if we look only at the part of $\dtr$ that lies above the diagonal line $\{ \rho=s\}$, then the maximal height in this subregion is bounded by $u$. 
\begin{lemma}\label{height bound}
Uniformly in the set of $r,t$ values lying in $\{ (r,t) : 0 \le r \le t\}$, we have that for any point $(\rho',s')\in\dtr \subset \R^+_{\rho'} \times \R^+_{s'}$,
\begin{enumerate}
\item
If $r \le t/3$, then
\begin{equation*}
|\dtr\cap \{(\rho',s') : \rho = \rho'\} | \le \min\{2\rho,2r\}
\end{equation*}
\item
If $t\ge r\ge t/3$, then
\begin{equation*}
|\{ s'\ge \rho'\ge0\} \cap \dtr\cap \{(\rho',s') : \rho = \rho'\} | \le u
\end{equation*}
\end{enumerate}
where $| \cdot|$ denotes the length. This implies that in either case, the height is bounded by $2\min(r,u)$. 
\end{lemma}
\begin{proof} We split the proof into two cases. 
\begin{enumerate}
\item
Let $r \le t/3$; then for each $\rho$, the maximal vertical length within $\dtr$ is $2r$ and occurs when $r \le \rho \le \f{t-r}2$; by symmetry, this length, $2r$, is maximal. When $0 \le \rho \le r$, the maximal vertical length of $\dtr$ is $2\rho$, which implies that this value of this length is sharp if and only if $0 \le \rho \le r.$

\item
Let $r \ge t/3$; then for each $\rho$, the maximal vertical length within $\dtr \cap \{ s \ge \rho \}$ is $t-r$ and occurs when $\f{t-r}2 \le \rho \le r$ and by symmetry once more, this length, $t-r$, is maximal. Furthermore, in a manner precisely analogous to the $r\le t/3$ case, we once more have that when $0 \le \rho \le \f{t-r}2$, the bound $2\rho$ is sharp if and only if $\rho$ lies in this small region.
\end{enumerate}
\end{proof}

\begin{proposition}[Application of the $r^\gamma$ estimate] \label{rp gain}
Let $\p$ solve 
\begin{equation}
\label{eqn}P\phi = \sum_{i=1}^{M} c_{i}(t,x) \phi^{2}\pa_{(i)}\p, \quad \pa_{(i)}\in \{ \paa \}_{\x \in \{0,\dots,3\}}, \quad c_{i}\in S^{Z}(1)
\end{equation}
for some $M\in\N$. 
Assume the hypotheses on $\gamma$ in \cref{the rp est} and also \cref{coeff.assu}. 
If $$\|\pa\pmn\|_{LE(T_1,T_2)} + \|\pa^2\pmn\|_{LE(T_1,T_2)} <\iy$$
for a sufficiently large $n$ and 
$\xg > 1/2$,
then
\begin{equation}\label{int.est}
|\jr (\p_3)_{\le m}| \ls \ju^{1/2 -\gamma/2}, \quad u > 1
\end{equation}
\begin{equation}
|(\p_3)_{\le m} | \ls r^{-(1 + \xg)/2}, \quad u < -1
\end{equation}
\end{proposition}

\begin{proof} We shall take $\xg>1/2$. 
Let $u > 1$. 
Let
$$H_3(s,\rho) := \sum_{k=0}^2 \| \xO^k (P\phi)_{\le m}(s,\rho \xo)\|_\ltst.$$ 
Recall \cref{eqn}, for any dyadic number $R$ (so that $\co$ can be either close to the origin or to the light cone). By \cref{the rp est} and Cauchy-Schwarz we have
\begin{align*}
\int_{\co}\rho H_{3}dA 
	&\ls \left( \int_{\co} \rho^{3-\xg}\|(\pmn)^{2}(\pa\pmn)^{2}\|_{\lt(S^{2})} \, dsd\rho\right)^{1/2}\\
	&\ls \left( \int_{\co}  \rho \jrho^{-\xg} \f1{\js^{2}}  \, dsd\rho\right)^{1/2}
\end{align*}The first line follows by \cref{the rp est}. The second line follows by \cref{u/v decay} (and the assumption of finite $LE^{1}$ norm), and \cref{derbound}. More precisely, we have
$$\rho^{3-\xg} \left( \f\ju{\jv^{2}} \right)^{2} \f1{\mu^{2}}$$
and we observe that $\mu\inv \sim \jv(\jr\ju)\inv$, which gives the integrand above.

We now split the collection of $R$ into those that are $\ll u$ and those otherwise. For the former set we obtain, by \cref{height bound},
$$\left( \sum_{R: R\ll u}\int_{\co}  \rho \jrho^{-\xg} \f1{\js^{2}}  \, dsd\rho\right)^{1/2}
	\ls \left( \sum_{R : R\ll u} R^{2-\xg} \ju\inv \right)^{1/2} \ls \ju^{1/2 - \xg/2}.$$
For the latter set, notice that for $\xg>1/2$ close to 1/2 we have $\jrho^{1-\xg} \le \js^{1-\xg}$, and we integrate over the remaining subset of $\dtr$ (call it $\dtr'$) to obtain
$$\left( \int_{\dtr'}  \rho \jrho^{-\xg} \f1{\js^{2}}  \, dsd\rho\right)^{1/2}
	\ls \left( \int_{\dtr'} \f1{\js^{1+\xg}} dsd\rho \right)^{1/2} \ls \ju^{1/2-\xg/2}.$$The proof for $u<-1$ is similar.
\end{proof}

\begin{remark}\label{rem:small gam}
We shall take $\xg= \f12+$ for \cref{eqn}. More precisely, given a fixed $\xs > \f14$ for \cref{eqn}, we pick $\xg$ such that the hypothesis $\xg<2\xs$ from the above theorem will be satisfied. 
\end{remark}

\section{The iteration in $\{ r > t +1 \}$}\label{sec:ext}

In this section we prove the pointwise decay rate stated in the main theorem in the region $\{r>t+1\}$.

\begin{theorem}
If $r>t+1$, then 
\begin{equation}\label{des.bd}
\pm \ls \jr\inv \ju^{-\min( 1 + \xs , \calT+\calN - 2  )}.
\end{equation}\end{theorem}

\begin{proof} %
	We shall assume only fairly weak bounds on the nonlinear term $H_3$.  The model we use for $H_3$ in the first part of the proof is $\bar\pa \p\pa\p$; all our nonlinearities decay at least as fast as this nonlinearity. In this sense, the first part of this proof will be catered to the linear part of the equation, which is the part that produces the $\ju^{-(1 + \xs)}$ bound. In the second part of our proof, we will prove the full decay rate for $\phi_3$ by using the full decay rate for $H_3$. 

Our initial bounds combined with \cref{tang.prop} can in $\{ r > t+1\}$ be written as
\begin{equation}\label{inbd}
\pmn\ls \frac{\ju^{1/2}}{\la r\ra}, \quad \pa\pmn \ls \frac{1}{\la r\ra \la u\ra^{1/2}},  \quad \overline{\partial} \phi_{\le m+n} \ls  \f{\ju^{1/2}}{\jr^2}.
\end{equation}
Since $\la u\ra\leq\la r\ra$, this can be weakened to 
\begin{equation}\label{1stbd}
\pmn\ls \frac{1}{\la r\ra^{1/2}}, \quad \pa\pmn \ls \frac{1}{\la r\ra^{1/2} \la u\ra}, \quad \overline{\pa} \phi_{\le m+n} \ls  \f{1}{\jr^{3/2}}.
\end{equation}
	Recall the decomposition \eqref{decomp}, and let 
$$H_k = \sum_{i=0}^2 \|\Omega^i (G_k)_{\leq m+n} (t, r \cdot)\|_{L^2(\S^2)}, \qquad k \in \{1,2,3\}.$$
	
	By \eqref{1stbd} (and \eqref{inbd} for $H_3$):
$$H_1 \ls \frac{1}{\la r\ra^{5/2+\sigma}}, \quad \pa_t H_2 
\ls \frac{1}{\la r\ra^{3/2+\sigma}\la u\ra}, \quad H_3 
\ls \frac{1}{\jr^{2+\kappa} \ju^{1 - \kappa}}, \quad \kappa \in (0,1).$$
For instance, for the nonlinearity of the form $\phi^{2}\pa\p$, by \cref{rp gain} and \cref{derbound} we have $$H_{3}(t,r)\ls \jv^{-2}\ju^{-\f{3\xg-1}2}$$ and this is bounded by $\jr^{-(2+\xk)}\ju^{-(1-\xk)}$ because $\xg$ is large enough.

	By \eqref{Bd2} with $\alpha=5/2+\sigma$, $\beta = 0$, and $\eta=0$, we obtain
\begin{equation}
\label{ineq1}(\p_1)_{\leq m+n} \ls \jr^{-1/2-\sigma} 
\end{equation}
which gains a factor of $\la r\ra^{-\sigma}$ compared to \eqref{1stbd}. Similarly \eqref{Bd2} with $\alpha=3/2+\sigma$, $\beta = 0$, and $\eta=1$ yields
\begin{equation}
\label{ineq2}(\p_2)_{\leq m+n} \ls \jr^{-1/2-\sigma}
\end{equation}
Next, \eqref{Bd2} with $\alpha=2+\xs$, $\beta = 0$, and $\eta=1/2$ yields
\begin{equation}
\label{ineq3}(\p_3)_{\leq m+n} \ls \jr^{-1/2-\sigma} 
\end{equation}
		\cref{ineq1,ineq2,ineq3} combined with \cref{derbound,tang.prop}, give the following improved bounds (by a factor of $\la r\ra^{-\sigma}$)
\begin{equation}\label{2ndbd}
  |\p_{\le m+n}| \ls \frac{1}{\jr^{1/2+\sigma}}, \quad |\pa\p_{\le m+n}| \ls \frac{1}{\jr^{1/2+\sigma} \la u\ra}, \quad \overline{\partial} \phi_{\le m+n} \ls \f{1}{\jr^{3/2+\sigma}}.
\end{equation}

We now repeat the iteration, replacing $\alpha$ by $\alpha+\sigma$ and applying \eqref{Bd2}. The process stops after some $k = k(\xs)$ steps, when \eqref{Bd2}, combined with \cref{derbound} and \eqref{tang.prop}, yield
\begin{equation}\label{3rdbd}
  |\p_{\le m+n}| \ls \frac{1}{\la r\ra}, \quad |\pa\p_{\le m+n}| \ls \frac{1}{\la r\ra \la u\ra}, \quad \overline{\partial} \phi_{\le m+n} \ls \f{1 }{\jr^2}.
\end{equation}

We now use only \eqref{Bd1} for $\p_1$ and $\p_3$, and \eqref{Bd1der} for $\p_2$. Note that \eqref{3rdbd} implies
\[
H_1 \ls \frac{1}{\la r\ra^{3+\sigma}}, \quad 
H_2 	\ls \frac{1}{\la r\ra^{2+\sigma}}, \quad 
H_3 	\ls \frac{1}{\la r\ra^3\la u\ra} .
\]For instance, for the nonlinearity of the form $\p^{2}\pa\p$, we the information in \cref{3rdbd} yields $H_{3}\ls \jr^{-3}\ju\inv$ exactly.

By \eqref{Bd1} with $\alpha=2+\sigma$, $\beta = 1$, and $\eta=0$, we obtain
\begin{equation}
\label{ineq4}(\p_1)_{\leq m+n} \ls \jr^{-1} \la u\ra^{-\sigma} 
\end{equation}
Similarly \eqref{Bd1der} with $\alpha=2+\sigma$, and $\eta=0$ yields
\begin{equation}
\label{ineq5}(\p_2)_{\leq m+n} \ls \jr^{-1} \la u\ra^{-\sigma} 
\end{equation}
Finally, \eqref{Bd1} with $\alpha=5/2$, $\beta = 1/2$, and $\eta=\sigma$ yields
\begin{equation}
\label{ineq6}(\p_3)_{\leq m+n} \ls \jr^{-1} \la u\ra^{-\xs} 
\end{equation}
	The bounds \cref{ineq4,ineq5,ineq6} combined with \cref{derbound,tang.prop} give the following improved bounds (by a factor of $\la u\ra^{-\sigma}$)
\begin{equation}\label{4thbd}
\pmn\ls \frac{1}{\la r\ra\ju^\sigma}, \quad \pa\pmn \ls \frac{1}{\la r\ra \la u\ra^{1+\sigma}}, \quad \overline{\pa} \pmn \ls \f1{\jr^2\ju^\sigma}.
\end{equation}
 
	We repeat this iteration and we can continue improving the decay rates of $\p_1$ and $\p_2$ to
	$$(\p_1)_{\le m+n}, \ \  (\p_2)_{\le m+n} \ls r^{-1} \ju^{-1}.$$
On the other hand, by \cref{4thbd} and the assumptions on C, it is clear that
$$H_3 \ls \frac{1}{\la r\ra^3 \la u\ra^{1+2\sigma}}$$
 and \eqref{Bd1} now yields
\begin{equation}
\label{yields}(\p_3)_{\le m+n} \ls r^{-1} \la u\ra^{-1}.
\end{equation}
One more iterate produces 
\begin{equation}
\label{light}(\p_1)_{\le m+n}, \ \  (\p_2)_{\le m+n} \ls r^{-1} \ju^{-1-\xs}.
\end{equation}
By the previous bounds we now have, for the \textit{original} value of $\xs$ from \cref{thm:main},
$$ H_1 \ls \f1{\jr^{3 + \xs}\ju^{1+}}, \quad H_2 \ls \f1{\jr^{2 + \xs}\ju^{1+}}, \quad H_3 \ls \f1{\jr^\calN \ju^{1+}}.$$
Using \cref{Bd1,Bd1der} now completes the proof for $\p_1,\p_2$. The proof for $\p_3$ is also complete if $\calT = 0$: an application of \cref{Bd1} with the above bound for $H_{3}$ shows that
$$(\p_{3})_{\le m+n}\ls r\inv\ju^{-(\calN - 2)}.$$

	We now consider the case when $\calT$ is nonzero. \textit{We note that the strategy, both here and in \cref{sec:int}, will always be to %
	use \cref{tang.prop} first for every tangential derivative in $\calC(\pa^{2}\p,\pa\p,\p)$, and then to use \cref{conversion}.}
Note that $\mu\inv = \min(\jr,\ju)\inv \sim \jv/(\jr\ju)$, and let $J:= \sum_{j=1}^\calN j_i$. 
By \cref{tang.prop} 
\begin{align*}
H_3 &\ls ( \ju\jv\inv)^{\calT} \mu^{- J} (\p_1)_{\le m+n} \dots (\p_\calN)_{\le m+n}\\
&\ls \jv^{-(\calT-J+\calN)} \ju^{-(J - \calT+\calN)} \jr^{-J} &\textrm{by }\cref{yields}
\end{align*}
where $\p_j$ is the $j$-th function in the nonlinearity.\footnote{The reader can verify that $\calT \le J$, so that near the light cone, we have decay in the $u$ variable of at least $\ju^{-\calN} \le \ju^{-3}$.}
	Since $u < -1$, we have $r\sim v$ and thus
\begin{align*}
H_{3} &\ls\jv^{-(\calT+\calN)} \ju^{-(J - \calT+\calN)}.
\end{align*}
	Observe that because $J-\calT\ge0$ and $\calN\ge3$, the value of $J-\calT+\calN$ lies in $(1,\iy)$, and so
by \cref{conversion}, we obtain
$$(\p_3)_{\le m+n} \ls r\inv \ju^{-(\calT+\calN-2)}.$$
In view of \cref{light}, this completes the proof because we have obtained
$$\pm\ls r\inv\ju^{-\min(1+\xs,\calT+\calN-2)}.$$
\end{proof}

\section{The iteration in $\{ r <t-1\}$} \label{sec:int}
\subsection{Converting $r$ decay to $t$ decay}\label{ss:inwardpropagation}

The pointwise decay rates for $\pmn$ obtained from the estimates for the fundamental solution are insufficient, but we show below that if the $\jr\inv$ decay from the fundamental solution is converted into $\jt\inv$, then the iteration does work.  

We will sometimes use the notation $C_{T_1}^{T_2} := [T_1,T_2]\times \{x : r\le t\}$.

\begin{lemma} Suppose that $\p$ satisfies the bound proved in \cref{led.for.vf}. Then for all $0 \leq  T_1 \leq T_2$, we have
\begin{align}\label{wwled est}
\begin{split}
\| \de \phi_{\leq m}&\|_{\lt(C_{T_1}^{T_2})}
\ls \sum_{j=1}^2 \| \jr^{1/2} \de \pm(T_j)\|_{\lt}  + \|\jr (P\p)_{\leq m}\|_{L^2} + \|\pat\pm\|_{\lt}.
\end{split}
\end{align}

\end{lemma}

\begin{proof}
We demonstrate the case $m=0$ first for simplicity. We multiply the equation by $r\pa_r\phi + \phi$ and integrate by parts in $[T_1,T_2]\times\R^3$. There is a number $q'>0$ such that
\begin{align}\label{computation}
\begin{split}
 \int |\de &\phi|^2 + O(\jr^{-q'})  |\de \phi|^2 + O(\jr^{-1-q'}) |\pao\p|^2 + O( \jr^{-2-q'}) |\phi|^2  \,dxdt \\
  &\ls \sum_{j=1}^2 \int_{\R^3} O(\la r \ra) |\de \phi(T_j,x)|^2 +  O(\la r \ra^{-1}) |\phi(T_j,x)|^2 \, dx + \int |r(P\p)\pa_r\phi| + |(P\p)\p| \,dxdt \\
  &\ls \sum_{j=1}^2 \int_{\R^3} O(\la r \ra) |\de \phi(T_j,x)|^2 \, dx + \int |r(P\p)\pa_r\phi| + |(P\p)\p| \,dxdt \\
\end{split}
\end{align}with the last statement following by a version of Hardy's inequality. 
	Next, by Cauchy-Schwarz and Hardy's inequality we can bound all the terms involving $P\p$ by 
	$$\f1\eps \|rP\p\|_\ltlt^2 + \eps\|\pa_r \phi\|_\ltlt^2.$$ 	
	By using the positivity of $q'$ on the left-hand side of \cref{computation} for large $|x|$ values, we can then obtain
\begin{equation}\label{sled m=0}
\begin{split}
\| \de \phi&\|_{\lt[T_1,T_2]\lt} \ls \sum_{j=1}^2 \| \jr^{1/2} \de \phi(T_j)\|_{L^2} +  \|\jr P\p\|_{L^2[T_1,T_2]L^2} .
\end{split}
\end{equation}
\cref{sled m=0} implies \eqref{wwled est} for $m=0$. 

\item (The higher multiindex case)
We now prove \eqref{sled m=0} but for $\phi_J, J\ne \vec 0$. 
We have
\begin{align*}
P \phi_J 	&= (P\p)_{J} + O(\jr^{-1-q'}) \de \phi_{\leq |J|} + O(\jr^{-2-q'}) \phi_{\leq |J|-1}.
\end{align*}
We multiply this by $r\pa_r \phi_J + \phi_J$. Then we integrate in $[T_1,T_2] \times \R^3$. The rest of the proof is then similar.
\end{proof}

\begin{lemma}\label{lem:aux}
Assume that $\p$ satisfies \cref{def:SILED}. Then 
$$\|\pm\|_{LE^1(\inte)} \ls T\inv \|\jr \pmn\|_{LE^1(\inte)} + \|(P\phi)_{\le m+n}\|_{LE^*(\inte)}.$$
\end{lemma}
\begin{proof}
Fix a dyadic number $T$. 
	Recall the SILED hypothesis, which we apply to dyadic time intervals $[T,2T]$, and recall \cref{led.for.vf} which states that we have SILED for vector fields. 
	We may assume that $\phi$ is supported in $\inte$ because we can control $[P,\chi]$, where $\chi$ is a purely spatial cutoff localised to the interior region $\{ r < 3t/4\}$, in the $LE^{*}$ norm. Since we assume SILED holds, we need not perform any cutoffs in the time variable. 

	Let $m\geq 0$. Let $\xg_{(T,x)}(t')$ denote an integral curve of $S$, parametrized by unit speed, such that $t'=0$ corresponds to the point $(T,x)$. 
	By the fundamental theorem of calculus and Cauchy-Schwarz, we have
\begin{equation}
|\de \pm(T,x)|^2
\ls \f1T \int_0^{T} 
|(\de \pm)(\gamma_{(T,x)}(t'))|^2 + |(S\de \pm)(\gamma_{(T,x)}(t'))|^2\,dt'
\label{vfromsv}\end{equation}
	A similar bound holds for $t = 2T$.
	Thus, after we integrate in $x$, we control the energy terms by 
	$$T^{-1/2}\|\pa \pmn\|_{L^{2}_{t,x}}.$$
	By using the fact that $\pat = t\inv(S - r\pa_{r})$,
	$$\|\pat \pm\|_{LE} \ls T\inv \|(S\pm, r\pa_{r}\pm)\|_{LE} \ls T^{-1/2}\|\jr\pmn\|_{\leo}.$$
	By \cref{led.for.vf} we now control the $\leo$ norm by 
		$$T^{-1/2}\|\pa \pmn\|_{L^{2}_{t,x}} + T^{-1/2}\|\jr\pmn\|_{\leo} + \|(P\p)_{\le m+n}\|_{\les}$$
	and we now use \cref{wwled est} to control the first term. By the fundamental theorem of calculus and Cauchy-Schwarz we have 
$$\|\jr^{1/2}\pa\p_{\le k}(T)\|_{\lt} \ls T^{-1/4}\|\jr^{1/4}\p_{\le k}\|_{\lt} + T^{-1/2}\|\jr\p_{\le k+1}\|_{\leo}$$
and similarly for the $t=2T$ energy norm. We decompose
$$\|\jr^{1/4}\de\p_{\le k}\|_{\lt} = \sum_{R < T} \|R^{1/4}\de\p_{\le k}\|_{\lt(r\sim R)}$$ and note that for all large $R$, 
$$\|R^{1/4}\de\p_{\le k}\|_{\lt(r\sim R)}\ls \|\jr\p_{\le k+n}\|_{\leo}$$
while for all sufficiently small $R$, we may absorb this to the left-hand side. 

On the other hand, 
\begin{align*}
\|\pat\p_{\le k}\|_{\lt} &\ls T\inv \|\p_{\le k+1}\|_{\lt} \\
	&\ls T^{-1/2} \sum_{R} \f{R^{1/2}}{T^{1/2}} \|\f{\p_{\le k+1}}{R^{1/2}}\|_{\lt(A_{R})}\\
	&\ls T^{-1/2} \sum_{R} \f{R^{1/2}}{T^{1/2}}\sup_{R} \|\f{\p_{\le k+1}}{R^{1/2}}\|_{\lt(A_{R})}\\
	&\ls T^{-1/2} \sum_{R} \f{R^{1/2}}{T^{1/2}} \|\p_{\le k+1}\|_{LE}\\
	&\ls T^{-1/2}  \|\p_{\le k+1}\|_{LE}.
\end{align*}
This concludes the proof.
\end{proof}

The next proposition uses \cref{lem:aux} to obtain better pointwise decay for the solution and its vector fields in the region $\{ r < t/2\}$. 

\begin{proposition}\label{convrsn}
Let $\phi$ solve \eqref{eq:problem}. Let $\xd > 0$. Assume that 
\begin{equation}\label{r bds}
\phi_{\le M}|_{r \le 3t/4} \ls \jr^{-1}\ju^{1/2 - q}, \quad \phi_{\le M}|_{r \le 3t/4} \ls \jt\inv \ju^{1/2 - q+\xd}, \quad q \geq \delta
\end{equation}
for an $M$ that is sufficiently larger than $m$.

If $q\geq \delta$ and $\xd \leq 1$, then
	$$ \pm|_\inte\ls  \jt\inv \ju^{1/2 - q} .$$
\end{proposition}

\begin{proof}
By \cref{derbound} and \cref{r bds},
$$ T\inv \|\jr \pmn\|_{LE^1} \ls T\inv \|\pmn\|_{LE} \ls T^{ - q}.$$

 Fix $n \geq 1$. By \cref{r bds}, and using \cref{2ndDeBd'} to bound each second-order derivative in the nonlinearity $P\phi$, we now have
 $$\|\jr^{1/2} (P\phi )_{\le m+n}\|_{L^2_{x,t}([T,2T]\times \{ r < 3t/4 \})} \ls T^{-q}.$$
That this is a weak bound for all nonlinearities other than $\phi^{2}\pa \p$ is not difficult to see. 

In what follows we prove that for the nonlinearity $\p^{2}\pa\p$, we obtain this bound.
	For $A_{R = 1}$, we use the $\jt\inv\ju^{1/2-q+\xd} \sim T^{-1/2-q+\xd}$ bound in \cref{r bds} to obtain 
 $$\|(\p^{2}\pa\p)_{\le m+n}\|_{L^2[T,2T]\lt(A_{R = 1} )} \ls T^{-1-3q+3\delta}.$$
 	 For $A_R, R>1$, we use \cref{derbound}, and we aim to obtain a power of $R^{-1.5}$ for the $R$ variable, so that after taking the volume element into account, the radial component of the norm will integrate to 1:
\begin{align*}
\|\jr^{1/2} (\p^2\pa\p)_{\le m+n}\|_{L^2[T,2T]\lt(A_R)} 
  &\ls \|R^{-1/2}(\pmn)^{3}\|\\
  &\ls \|R^{-1/2} (R\inv T^{1/2-q}) (T^{-1/2-q+\xd})^{2}\| \\
  &\ls T^{1 - q} T^{-1 - 2q + 2\xd} = T^{-3q + 2\delta}\\
  &\ls T^{-q} \quad \text{if and only if } q \ge \delta.
\end{align*}
 
Therefore \cref{lem:aux} implies
$$\|\pm\|_{LE^1(\inte)} \ls T^{- q}$$
and the conclusion now follows by \cref{DyadLclsd}. We note that the bound in $A_{R=1}$ is no worse than the bound $T^{-q}$ so long as $\xd \leq 1$; this is the only time we use the hypothesis $\xd\le1$. 
\end{proof}

\begin{remark}
In the present article we shall let
$$\xd := \min(\xg/2, \xs).$$
By \cref{conversion} applied to the linear components of the equation (i.e. the coefficients of $P-\Box$), we obtain
$$\jr(\p_{1})_{\le m+n}, \jr(\p_{2})_{\le m+n}\ls \ju^{1/2-\xs}.$$
By \cref{rp gain}, we obtain
$$\jr(\p_{3})_{\le m+n}\ls \ju^{1/2-\xg/2}.$$
These two estimates imply that 
$$\jr\pmn\ls \ju^{1/2-\delta}.$$
Thus the hypothesis $q \ge \delta$ in \cref{convrsn} holds throughout our iteration for nonlinearities containing terms of the form $S^{Z}(1)\p^{2}\pa\p$.
\end{remark}

\begin{remark}
We shall use \cref{convrsn} whenever the vector fields $\pmn$ satisfy a bound of the form $\jr\inv D$ where $D$ is some decaying function. Bounds of the form $\jr\inv D$ arise when applying the fundamental solution for $\Box$, which is our strategy in the iteration (see \cref{the iter}). \cref{convrsn} then turns this upper bound into the upper bound $\pmn\ls \jt\inv D$. 
\end{remark}

\subsection{The iteration}\label{the iter}

\begin{theorem}
If $r<t-1$, then 
\begin{equation}
\label{full}\pm \ls \jv\inv \ju^{-\min(1+\xs,\calT+\calN-2)}.
\end{equation}
Here $\xs$ denotes the original value of $\xs$ taken from \cref{thm:main}.
\end{theorem}

\begin{proof}Just as in the proof of \cref{des.bd}, in the first part of this proof, we shall assume only fairly weak bounds on the nonlinear term $H_3$.  The model we use for $H_3$ in the first part of the proof is $\bar\pa \p\pa\p$; all our nonlinearities in \cref{eq:problem} decay at least as fast as this nonlinearity. %
The idea is that we wish \begin{enumerate}
\item
first to establish a decay rate of $\pmn\ls\jv\inv\ju\inv$, and 
\item
second to establish the full decay rate \cref{full}. 
\end{enumerate}

\noindent \textit{Step 1.} By \cref{tang.prop} and our initial decay estimates we have, in $\{ r<t/2\}$,
\begin{equation}\label{inbd1}
\pmn\ls \f{\ju^{1/2}}{\la t\ra}, \quad \pa\pmn \ls \f{1}{\la r\ra \la u\ra^{1/2}},  \quad \overline{\partial} \phi_{\le m+n} \ls \f{\ju^{1/2}}{\jr\jt}.
\end{equation}
	These bounds \cref{inbd1} imply
$$H_1 \ls \frac{\ju^{1/2}}{\jr^{2+\xs}\jt}, \quad \pa_t H_2 \ls \frac{1}{\la r\ra^{1+\sigma}\jt\la u\ra^{1/2}}, \quad H_3 \ls \frac{1}{\la r\ra^2\la t\ra}.$$
For the equation \cref{eqn}, the result of \cref{rp gain} implies that $H_{3}$ for that equation satisfies
$$H_{3} \ls \mu\inv \cdot (\pmn)^{3} \ls \f1{\jt^{2}\jr\ju^{\xl}}, \quad \xl = \f14 + 3\eps$$
if we take $\xg = \f12 + \eps$ (recall \cref{rem:small gam}). Thus $H_{3} \ls \jr^{-2}\jt\inv$. 
	On the other hand, for other nonlinearities, direct computation shows that $H_{3}\ls \jr^{-2}\jt\inv$. 

	By \eqref{Bd1} with $\x=2+\sigma$, $\beta = 1$, and $\eta=-1/2$, we obtain 
\[
(\p_1)_{\leq m+n} \ls \jr^{-1} \ju^{1/2-\sigma} %
\]
	Similarly \eqref{Bd1} with $\alpha=2+\sigma$, $\beta = 0$, and $\eta=1/2$ yields
\[
(\p_2)_{\leq m+n} \lesssim \jr^{-1} \ju^{1/2-\sigma} %
\]
	Finally, \cref{Bd1} with $\x=2+\sigma$, $\beta = 1-\sigma$, and $\eta=0$ yields
$$(\p_3)_{\leq m+n} \lesssim \jr^{-1} \ju^{1/2-\sigma}.$$
Thus $\pmn\ls \jr^{-1} \ju^{1/2-\xs},$
which represents our first improvement over the initial pointwise decay rate for $\phi$.

By \cref{convrsn} with $q=-1/2+\eta$, and \cref{derbound,tang.prop} we obtain the following improvement, by a factor of $\ju^{-\sigma}$, over \cref{inbd1}:
\begin{equation}\label{inbd2}
\pmn \ls \frac{\ju^{1/2-\sigma}}{\la t\ra}, \quad \pa\pmn \ls \frac{1}{\la r\ra \la u\ra^{1/2+\sigma}},  \quad \overline{\partial} \phi_{\le m+n} \ls  \f{\ju^{1/2-\xs}}{\jr\jt}.
\end{equation}

\noindent \textit{Step 2.} We iterate, replacing $\eta$ by $\eta+\sigma$, applying \eqref{Bd2} and then using \cref{convrsn} to turn the $\jr\inv$ factor in the upper bound for $\pmn$ into a $\jt\inv$ factor. The process stops after $\lfloor(2\xs)\inv\rfloor$ steps, when \eqref{Bd1}, combined with \cref{convrsn,derbound,tang.prop}, yield
\begin{equation}\label{3rdbd'}
\pmn\ls \frac{1}{\la t\ra}, \quad \pa\pmn \ls \frac{1}{\la r\ra \la u\ra}, \quad \overline{\partial} \pmn \ls \f{1}{\jr\jt}.
\end{equation}

\noindent \textit{Step 3.} We use \cref{Bd1der} for $\phi_2$, and the iteration process follows the same pattern as in \cref{sec:ext}, but with \cref{convrsn} to turn upper bounds for $\pmn$ of the form $\jr\inv (\dots)$ into upper bounds of the form $\jv\inv(\dots)$. We obtain 
$$(\p_{1})_{\le m+n}, (\p_{2})_{\le m+n} \ls \jv\inv \ju^{-1-\xs}.$$
	This completes the proof of Step (1) above; that is, we have 
\begin{equation}\label{to-beat}
\pmn\ls\jv\inv\ju\inv.
\end{equation} 

\noindent \textit{Step 4.} We now prove the full decay rate \cref{full}. All that remains is to prove the bound
$$(\p_{3})_{\le m+n}\ls\jr\inv\ju^{-(\calT+\calN-2)}$$
and an application of \cref{convrsn} then completes the proof. 

Recall that $\mu\inv = \min(\jr,\ju)\inv \sim \jv/(\jr\ju)$, and let $J:= \sum_{j=1}^\calN j_i$. 
Then
\begin{align*}
H_3 &\ls ( \ju\jv\inv)^{\calT} \mu^{- J} (\p_1)_{\le m+n} \dots (\p_\calN)_{\le m+n} &\textrm{by }\cref{tang.prop}\\
	&\ls \jv^{-(\calT-J+\calN)} \ju^{-(J - \calT+\calN)} \jr^{-J} &\textrm{by }\cref{yields}
\end{align*}
Recall the requirement that $\xb\ge0$ in the hypothesis of \cref{conversion}. Thus we want the exponent of $\ju$ in the display immediately above to be non-negative. To this end, we replace $J$ by $J' := \min(J,3)$, that is to say we apply \cref{derbound} at most three times; thus we now obtain 
\begin{equation}
H_{3}\ls \jv^{-(\calT-J' +\calN)} \ju^{-(J' - \calT+\calN)} \jr^{-J'}.
\end{equation}
	Then we apply \cref{Bd1}: there are two cases in that lemma, namely $\eta<1$ and $\eta>1$. The former gives a strictly worse bound (by a factor of $\ju$); we use this worse bound now to obtain, in the notation of that lemma,
\begin{equation}\label{update}
r(\p_{3})_{\le m+n}\ls 1/\ju^{\x + \xb+\eta-3} = 1/\ju^{2\calN + J'-3}.
\end{equation}
Note that 
$$2\calN+J'-3 > 1$$
so that, upon application of \cref{convrsn}, \cref{update} becomes a strictly better bound than \cref{to-beat}, namely
\begin{equation}
\label{btr}\pmn\ls \jv\inv\ju^{-\xk'}, \quad \xk' := \min(1+\xs, 2\calN + J'-3).
\end{equation}
	Next, by \cref{btr} we have
\begin{align*}
H_{3} &\ls \jr^{-J}\jv^{-(\calT - J+\calN)} \ju^{-(J - \calT + \xk'\calN)}\\
 &= \jr^{-J}\jv^{-(\calT - J+\calN)} \ju^{-\big(J - \calT +\min(\calN(1+\xs),  \calN(2\calN + J'-3) )\big)}\\
 &\le \jr^{-J}\jv^{-(\calT - J+\calN)} \ju^{-3-}
\end{align*}
and there is now no need to define a truncation $J'$ when attempting to apply \cref{conversion}. The final line follows because $J-\calT\ge0$ and $1+\xs, 2\calN + J'-3>1$ are both strictly larger than 1. 

	We apply \cref{Bd1} one final time, noting that $\eta>1$ in the notation of \cref{Bd1}, to obtain
$$r(\p_{3})_{\le m+n}\ls 1/\ju^{\x+\xb-2} = 1/\ju^{J + (\calT - J + \calN) - 2} = 1/\ju^{\calT+\calN-2}.$$
	Thus we conclude
$$\pmn \ls r\inv \ju^{-\min(1+\xs,\calT+\calN-2)}.$$
Applying \cref{convrsn} %
 now completes the proof, since we obtain
$$\pm \ls \jv\inv \ju^{-\min(1+\xs,\calT+ \calN-2)}.$$
\end{proof}

\section{Nonlinearities with a special structure \cref{eq:eqn,eq:eqn2}}\label{sec:final}

\begin{theorem}[Improved decay for equations \cref{eq:eqn,eq:eqn2} in spherical symmetry]
Let $\p(t,r)$ solve the equation \cref{eq:eqn} 
$$P\phi (t,r)= \sum_{i=1}^{M} c_{i}(t,r) \phi^{n}\pa_{(i)}\p, \quad \pa_{(i)}\in \pa \quad c_{i}\in S^{Z}_{\pa}(1), \quad n \ge2$$ or the equation \cref{eq:eqn2}.
Then
$$\pm\ls%
\jv\inv\ju^{-\min(1+\xs,n)}.$$
\end{theorem}

\begin{proof}
From previous work \cref{full,des.bd} we already have 
\begin{equation}
\label{usethis}
\pmn(t,r) \ls\jv\inv\ju\inv.
\end{equation}We only write down the proof for \cref{eq:eqn} because the proof for \cref{eq:eqn2} is identical.

We first do the constant $c_i$ case, since the general case then follows easily. We set the constant to equal 1/($n+1$), thus $P \p = \pa_{(i)}(\p^{n+1})$. 
\begin{itemize}
\item
In the case $\pa_{(i)} = \pa_r$, we note that the under the assumption of spherical symmetry (of the solution), and specialising to the Minkowski spacetime case for now (thus $P = \Box$), one has the equality 
\begin{equation}
\label{refhere}r\p(t,r) = \f12   \int_0^t \int_{|r-(t-s)|}^{r + (t-s)} \rho \pa_\rho(\p^{n+1}) dA .
\end{equation} 
and to upper-bound $r\p$ it suffices to bound the integral.
We set $n=2$ and note that the higher $n$ case follows similarly. Integrating by parts, we see that this equals 
$$\int_0^t \rho\p^3|_{\rho = |r-(t-s)|}^{r+(t-s)} ds - \int_\dtr\p^3 \,dA = I + II + III.$$
By \cref{usethis}
$$I\ls v^{-3} \int_0^t \f{r + t - s}{\la r + t - 2s\ra^3} ds \ls v^{-3} \cdot \f\jv\ju.$$
Also by \cref{usethis},
$$II\ls u^{-3} \int_0^t \f{|s-u|}{\la 2s - u\ra^3} ds \ls u^{-3} \cdot \ju\inv.$$
Finally, by \cref{usethis} note that $\p^3 \ls \jv^{-3}\ju^{-3}$ so in the notation of \cref{conversion}, the integrand $\p^{3}$ satisfies $\eta = 3, \xb + \x = 3$. Thus by \cref{conversion} we have
$$\int_\dtr\p^3dA\ls\ju^{-2}.$$
We conclude that 
$$r\p(t,r)\ls \ju^{-2}.$$
For general $n\ge2$, we have
$$\int_\dtr\p^{n+1}dA\ls \ju^{-n}$$
and 
$$r\p(t,r)\ls\ju^{-n}.$$

For the non-constant $c_i$ case, we note that the assumption $\pa c \in S^Z(\jr\inv)$ together with \cref{usethis} implies the same overall $\ju^{-2}$ bound. 

In the case $P \ne \Box$, the equation \cref{refhere} is modified into an equation that includes  integrals on the right hand side that depend on the coefficients of $P-\Box$:
$$r\p(t,r)= \f12\int_{\dtr} \rho\pa_{\rho}(\p^{n+1}) \,dA + \f12 \int_{\dtr} \rho V(s,\rho)\p  \,dA + \dots$$
where $\dots$ denotes other terms from $P-\Box$. 
The arguments for these coefficients of $P - \Box$ were shown in \cref{sec:ext,sec:int} already.
	
	Applying \cref{convrsn} completes the proof because the $r$ factor on the left hand side of \cref{refhere} is then converted into a $t+r$ factor, giving decay inside the interior region $\{ r<t/2\}$. 
\item
In the case $\pa_{(i)} = \pat$, we provide a proof that also works outside of spherical symmetry for %
the equation $P\p(t,x) = c\p^n\pat\p$. %

We once again first assume constant $c_i$. Again we first look at the case $n=2$. The idea is that in the interior $\{ r < t/2\}$ we may use the relation $\pat = t\inv(S -r\pa_r)$ and \cref{derbound} to bound $\pa_r$ to obtain extra $\js\inv$ decay when applying \cref{conversion}; this implies extra $\ju\inv$ decay. Rigorously, this means that we let $1 = \chi_\text{int}+\chi_\text{con}$ with $\chi_\text{int}$ localised to the interior. Then 
$$\pat(\p^3) = 3\chi_\text{int} \p^2\pat\p + \chi_\text{con}\pat(\p^3) =3 \chi_\text{int} \p^2\pat\p + \pat(\chi_\text{con}\p^3) - \pat \chi_\text{con} \p^3 = A+B+C.$$
Let $\Box\p_4 = A, \Box w = \chi_\text{con}\p^3, \Box \p_6 = C$. 
Note that 
$$\chi_\text{int}\p^2\pat\p \ls |\pn|^3 / \jt\inv$$
and applying \cref{conversion} as in previous sections gives the desired $\ju^{-2}$ bound.
On the other hand, we apply the procedure outlined in \cref{Minkdcyt} to obtain extra $\ju\inv$ decay for 
$$\Box w = \chi_\text{con}\p^3.$$
Thus $\jr \pat w \ls \ju^{-2}$. Since $\pat \chi_\text{con} \ls \jt\inv$, by applying \cref{conversion} we obtain $r\p_6 \ls \ju^{-2}$. In summary, by the triangle inequality we conclude that the solution $\p$ to the wave equation $\Box \p = \pat(\p^{3})$ with zero initial data obeys the pointwise bounds $\jr\p\ls\ju^{-2}$. Applying \cref{convrsn} completes the proof, since we now have $\jv\p\ls\ju^{-2}$.

The general $n$ case is very similar. The non-constant $c_i$ case follows by the assumption that for all $i$, $\pa c_{i}\in S^Z(\jr\inv)$. 
\end{itemize}
\end{proof}

\section*{Acknowledgements}
I would like to thank Siyuan Ma for a conversation.

\end{document}